\documentclass{amsart}
\usepackage{graphicx}
\usepackage{amsmath,amsthm}
\numberwithin{equation}{section}

\def\RR{\mathbb{R}}
\def\Q{\mathbf{Q}}

\def\N{\mathbb{N}}
\def\K{\mathbf{K}}
\def\U{\mathbf{U}}

\def\s{\mathbf{S}}
\def\P{\mathbf{P}}
\def\X{\mathbf{X}}
\def\U{\mathbf{U}}
\def\Si{\mathbf{\Sigma}}

\def\m{\mathcal{M}}
\def\I{{\rm I}}

\newtheorem{theorem}{Theorem}[section]

\newtheorem{proposition}[theorem]{Proposition}

\theoremstyle{definition}
\newtheorem{defn}[theorem]{Definition}

\theoremstyle{remark}
\newtheorem{remark}[theorem]{Remark}

\title[Optimal control and LMI-relaxations]{Nonlinear optimal control
via occupation measures and LMI-relaxations}

\author{Jean B.~Lasserre}
\address{LAAS-CNRS and Institute of Mathematics, 
University of Toulouse, France.}
\email{lasserre@laas.fr}

\author{Didier~Henrion}
\address{LAAS-CNRS, University of Toulouse, France and
Czech Technical University in Prague, Czech Republic.}
\email{henrion@laas.fr}

\author{Christophe~Prieur}
\address{LAAS-CNRS, University of Toulouse, France.}
\email{cprieur@laas.fr}

\author{Emmanuel~Tr\'elat}
\address{University of Orl\'eans, France}
\email{emmanuel.trelat@univ-orleans.fr}


\date{\today}


\begin{document}

\begin{abstract}
We consider the class of nonlinear optimal control problems (OCP) with
polynomial data, i.e., the differential equation, state and control
constraints and cost are all described by polynomials, and more generally for 
OCPs with smooth data.
In addition, state constraints as well as state and/or action
constraints are allowed. We provide a simple hierarchy of
LMI (linear matrix inequality)-relaxations whose optimal values form a nondecreasing sequence of
lower bounds on the optimal value. Under some convexity assumptions,
the sequence converges to the optimal value of the OCP. Preliminary
results show that good approximations are obtained with few moments.
\end{abstract}

\maketitle

\section{INTRODUCTION}

Solving a general nonlinear optimal control problem (OCP) 
is a difficult challenge, despite powerful theoretical tools are
available, e.g. the maximum principle and
Hamilton-Jacobi-Bellman (HJB) optimality equation.
The problem is even more difficult in the presence of state and/or control constraints.
State constraints are particularly difficult to handle, and the interested reader
is referred to Capuzzo-Dolcetta and Lions \cite{lions} and Soner \cite{soner}
for a detailed account of HJB theory in the case of state constraints.
There exist many numerical methods to compute the solution
of a given optimal control problem; for instance,
{\it multiple shooting} techniques which solve
two-point boundary value problems as described, e.g., in
\cite{sb02,Pesch}, or {\it direct
methods}, as, e.g., in \cite{sb,fl,gill}, which use, among others,
descent or gradient-like algorithms.
To deal with optimal control problems with state constraints, some
adapted versions of the maximum principle have been developed (see
\cite{Jacobson,Maurer}, and see \cite{Hartl} for a survey of this
theory), but happen to be very hard to implement in general.

On the other hand,  the OCP can be written as an infinite-dimensional linear program (LP) over two spaces of measures. This is called the {\it weak} formulation of the OCP in Vinter \cite{vinter} (stated in the more general context of differential inclusions).
The two unknown measures are the state-action {\it occupation measure} (o.m.) {\it up to} the final time $T$, and the state o.m. {\it at} time $T$.  The optimal value of the resulting LP always provides a
lower bound on the optimal value of  the OCP, and under some convexity assumptions, 
both values coincide; see Vinter \cite{vinter}  and Hernandez-Hernandez et al. \cite{hernandez} as well.

The dual of the original infinite dimensional LP has an interpretation
in terms of "subsolutions" of related HJB-like optimality conditions,
as for the unconstrained case.
The only difference with the unconstrained case is the underlying
function space involved, which directly incorporate the state
constraints. Namely, the functions are only defined on the state
constraint set .

An interesting feature of this LP approach with o.m.'s is that state
constraints, as well as
state and/or action constraints are all easy to handle; indeed they simply
translate into constraints on the supports of the unknown o.m.'s.
It thus provides an alternative to the use of maximum principles with
state constraints. However, although this LP approach is valid for any
OCP, solving the corresponding (infinite-dimensional) LP is difficult
in general; however, general LP approximation schemes based on grids
have been proposed in e.g. Hernandez and Lasserre \cite{lasserre}. 

This LP approach has also been used in the context of 
discrete-time Markov control processes, and is dual to Bellman's optimality principle.
For more details the interested reader is referred to Borkar \cite{borkar},
Hernandez-Lerma and Lasserre \cite{ohlbook,ohlbook2,ohllinear}
and many references therein. 
For some continuous-time stochastic control problems (e.g., modeled by
diffusions) and optimal stopping problems, the LP approach has also been used with 
success to prove existence of stationary 
optimal policies; see for instance Cho and Stockbridge \cite{stock0}, Helmes and Stockbridge \cite{helmes1}, Helmes et al. \cite{helmes}, Kurtz and Stockbridge \cite{kurtz}, and also
Fleming and Vermes \cite{fleming}. In some of these works, the
moment approach is also used to approximate the resulting infinite-dimensional LP.

\medskip

{\bf Contribution.} In this paper, we consider the particular class
of nonlinear OCP's with state and/or control constraints, for which all
data describing the problem (dynamics, state and control constraints)
are {\it polynomials}. The approach also extends to
the case of problems with {\it smooth} data and compact sets, because polynomials are dense
in the space of functions considered; this point of view is detailed in \S\ref{smooth}.
In this restricted polynomial framework, the LP approach has interesting
additional features that can be exploited for effective numerical
computation. Indeed, what makes this LP approach attractive is that 
for the class of OCPs considered:

$\bullet$ Only the {\it moments} of the o.m.'s appear in the LP formulation, so that we already end up with countably many variables, a significant progress. 

$\bullet$ Constraints on the support of the o.m.'s translate easily into 
either LP or SDP (Semi Definite Programming) {\it necessary}
constraints on their moments. Even more, for (semi-algebraic) compact supports,
relatively recent powerful results from real algebraic geometry make these constraints also {\it sufficient}.

$\bullet$ When truncating to finitely many moments, the resulting LP
or SDP's are solvable and their optimal values form a monotone
nondecreasing sequence (indexed by the number of moments considered)
of lower bounds on the optimal value of the LP (and thus of the OCP).

Therefore, based on the above observations, we propose an approximation of the optimal value of
the OCP via solving a hierarchy of SDPs (or linear matrix inequalities, LMIs)that provides a monotone nondecreasing sequence of lower bounds on the optimal value of the weak LP formulation of the OCP.
In adddition, under some compactness assumption of the state and control constraint sets,
the sequence of lower bounds is shown to converge to the optimal value of the LP, and thus the optimal value of the OCP when the former and latter are equal.

As such, it could be seen as a complement to the above shooting or
direct methods, and when the sequence of lower bounds converges to the
optimal value, a test of their efficiency.
Finally this approach can also be used to provide a {\it certificate} of unfeasibility. Indeed, if in the hierarchy of LMI-relaxations of the minimum time OCP, one is infeasible then the OCP itself is infeasible. It turns out that sometimes this certificate is provided at
an early stage in the hierarchy, i.e. with very few moments. This is illustrated on two simple examples.

In a pioneering paper, Dawson \cite{dawson} had suggested the use of {\it moments}
in the LP approach with o.m.'s,  but results on the $\K$-moment problem 
by Schm\"udgen \cite{schmudgen}  and Putinar \cite{putinar} were
either not available at that time. Later, Helmes and Stockbridge \cite{helmes1} and
Helmes, R\"{o}hl and Stockbridge
\cite{helmes} have used LP moment conditions for computing
some exit time moments in some diffusion model,
whereas for the same models, Lasserre and Prieto-Rumeau
\cite{lasserprieto} have shown that SDP moment conditions are superior in terms
of precision and number of moments to consider; in \cite{finance},
they have extended the moment approach for options pricing problems in
some mathematical finance models.
More recently, Lasserre, Prieur and Henrion \cite{cdc05} have used
the o.m. approach for minimum time OCP without state
constraint. Preliminary experimental results on
Brockett's integrator example, and the double integrator show fast
convergence with few moments.

\section{Occupation measures and the LP approach}
\label{general}
\subsection{Definition of the optimal control problem}
Let $n$ and $m$ be nonzero integers. Consider on $\RR^n$ the control system
\begin{equation}\label{system}
\dot{x}(t)=f(t,x(t),u(t)),
\end{equation}
where $f:[0,+\infty)\times \RR^n \times \RR^m\longrightarrow \RR^n$ is smooth,
and where the controls are bounded measurable functions, defined on intervals
$[0,T(\mathbf{u})]$ of $\RR^+$, and taking their values in a
\textit{compact} subset $\U$ of $\RR^m$.
Let $x_0\in \RR^n$, and let $\X$ and $\K$ be \textit{compact}
subsets of $\RR^n$.
For $T>0$, a control $u$ is said \textit{admissible} on $[0,T]$
whenever the solution $x(\cdot)$ of (\ref{system}), such that
$x(0)=x_0$, is well defined on $[0,T]$, and satisfies
\begin{equation}\label{cont1}
(x(t),u(t))\in \X\times\U\quad\textrm{a.e. on}\ [0,T],
\end{equation}
and
\begin{equation}\label{cont2}
x(T)\in \K.
\end{equation}
Denote by $\mathcal{U}_T$ the set of \textit{admissible controls} on
$[0,T]$.

For $\mathbf{u}\in \mathcal{U}_T$, the \textit{cost} of the associated
trajectory $x(\cdot)$ is defined by
\begin{equation}\label{cost}
J(0,T,x_0,\mathbf{u})=\displaystyle\int_0^T h(t,x(t),u(t))dt 
+ H(x(T)),
\end{equation}
where $h:[0,+\infty)\times\RR^n\times\RR^m \longrightarrow \RR$
and $H:\RR^n\rightarrow \RR$ are smooth functions.

Consider the \textit{optimal control problem} of determining a
trajectory solution of (\ref{system}, starting from $x(0)=x_0$,
satisfying the state and control constraints (\ref{cont1}), the
terminal constraint (\ref{cont2}), and minimizing the cost
(\ref{cost}). The final time $T$ may be fixed or not.

If the final time $T$ is fixed, we set
\begin{equation}\label{crit1}
J^*(0,T,x_0)\,:=\,\inf_{\mathbf{u}\in \mathcal{U}_T}
J(0,T,x_0,\mathbf{u}),
\end{equation}
and if $T$ is free, we set
\begin{equation}\label{crit1libre}
J^*(0,x_0)\,:=\,\inf_{T>0,\ \mathbf{u}\in \mathcal{U}_T}
J(0,T,x_0,\mathbf{u}),
\end{equation}
Note that a particular OCP is the minimal time problem
from $x_0$ to $\K$, by letting $h\equiv 1$, $H\equiv 0$. In this
particular case, the minimal time is the first {\it hitting time} of
the set $\K$.


It is possible to associate a stochastic or deterministic OCP with an abstract
infinite dimensional linear 
programming (LP) problem $\P$ together with its dual $\P^*$ (see for instance
Hern\'andez-lerma and Lasserre \cite{ohlbook} for discrete-time Markov
control problems, and Vinter \cite{vinter}, Hernandez et
al. \cite{hernandez} for deterministic optimal control problems, 
as well as many references therein). We next describe this LP approach
in the present context.

\subsection{Notations and definitions}
For a topological space $\mathcal{X}$, let 
$\m(\mathcal{X})$ be the 
Banach space of finite signed Borel measures on $\mathcal{X}$, equipped with the norm of total variation, and denote by $\m(\mathcal{X})_+$ its positive cone, that is, the space of finite measures 
on $\mathcal{X}$. Let $C(\mathcal{X})$ be the Banach space of bounded continuous functions on $\mathcal{X}$, equipped with the sup-norm. Notice that when $\mathcal{X}$ is compact Hausdorff, then
$\m(\mathcal{X})\simeq C(\mathcal{X})^*$, i.e., 
$\m(\mathcal{X})$ is the topological dual of $C(\mathcal{X})$.

Let $\Si:=[0,T]\times \X$, $\s:=\Si\times \U$, and let $C_1(\Si)$ be the Banach space of functions $\varphi\in C(\Si)$ that are continuously differentiable.
For ease of exposition we use the same notation 
$g$ (resp. $h$) for a polynomial  $g\in\RR[t,x]$ (resp. $h\in\RR[x]$) and its restriction to 
the compact set $\Si$ (resp. $\K$). 

Next, with $u\in \U$, let $A:\,C_1(\Si)\to C(\s)$ be the mapping
\begin{equation}
\label{mapping}
\varphi\mapsto\,A\varphi (t,x,u)\,:=\,\frac{\partial \varphi}{\partial
t}(t,x)+\langle f(t,x,u),\nabla_x \varphi(t,x)\rangle.
\end{equation}
Notice that $\partial\varphi/\partial t+\langle
\nabla_x\varphi,f\rangle\in C(\s)$ for every $\varphi\in C_1(\Si)$,
because both $\X$ and $\U$ are compact, and $f$ is understood as its
restriction to $\s$.

Next, let $\mathcal{L}:\,C_1(\Si)\to C(\s)\times C(\K)$ be the linear mapping
\begin{equation}
\label{map}
\varphi\,\mapsto\,\mathcal{L}\varphi \,:=\,(-A\varphi, \varphi_T),
\end{equation}
where $\varphi_T(x):=\varphi(T,x)$, for all $x\in \X$.
Obviously, $\mathcal{L}$ is continuous with respect to the strong
topologies of $C_1(\Si)$ and $C(\s)\times C(\K)$.

Both $(C(\s),\m(\s))$ and $(C(\K),\m(\K))$ form a {\it dual pair} of vector spaces, with duality brackets
\[\langle h,\mu\rangle =\int h\,d\mu,\qquad \forall \,(h,\mu)\in C(\s)\times \m(\s),\]
and
\[\langle g,\nu\rangle =\int g\,d\nu,\qquad \forall \,(g,\nu)\in C(\K)\times \m(\K) \]. 

Let $\mathcal{L}^*:\,M(\s)\times M(\K)\to C_1(\Sigma)^*$ 
be the adjoint mapping of $\mathcal{L}$, defined by
\begin{equation}
\label{map*}
\langle (\mu,\nu),\mathcal{L}\varphi\rangle\,=\,
\langle \mathcal{L}^*(\mu,\nu),\varphi\rangle,
\end{equation}
for all $((\mu,\nu),\varphi)\in
M(\s)\times M(\K)\times C_1(\Si)$.

\begin{remark}
\label{continuity}
\begin{itemize}
\item[(i)] The mapping $\mathcal{L}^*$ is continuous with respect to
the weak topologies $\sigma (\m(\s)\times \m(\K),C(\s)\times C(\K))$,
and $\sigma (C_1(\Si)^*,C_1(\Si))$.
\item[(ii)] Since the mapping $\mathcal{L}$ is continuous in the
  strong sense, it is also continuous with respect to
the weak topologies $\sigma (C_1(\Si),C_1(\Si)^*)$ and
$\sigma (C(\s)\times C(\K),\m(S)\times \m(\K))$.
\item[(iii)] In the case of a {\it free} terminal time $T\leq T_0$, it
  suffices to redefine
$\mathcal{L}:\,C_1(\Si)\to C(\s)\times C([0,T_0]\times \K)$ by
$\mathcal{L}\varphi :=(-A\varphi,\varphi)$. The operator
$\mathcal{L}^*:\,M(\s)\times M([0,T_0]\times \K)\to C_1(\Sigma)^*$ is
still defined by (\ref{map*}),
for all $((\mu,\nu),\varphi)\in
M(\s)\times M([0,T_0]\times\K)\times C_1(\Si)$.

For time-homogeneous free terminal time problems, one only needs functions $\varphi$
of $x$ only, and so $\Si=\s=\X\times \U$ and 
$\mathcal{L}:\,C_1(\Si)\to C(\s)\times C( \K)$.
\end{itemize}
\end{remark}

\subsection{Occupation measures and primal LP formulation}
Let $T>0$, and
let $\mathbf{u}=\{u(t),\:0 \leq t<T\}$ be a control such that the
solution of (\ref{system}), with $x(0)=x_0$, is well defined on $[0,T]$.
Define the probability measure $\nu^{\mathbf{u}}$ on $\RR^n$, and the measure
$\mu^{\mathbf{u}}$ on $[0,T]\times\RR^n\times \RR^m$, by
\begin{eqnarray}
\label{nu}
\nu^{\mathbf{u}}(D)&:=&\I_D\,[x(T)],\quad D\in\mathcal{B}_n,\\
\label{mu}
\mu^{\mathbf{u}}(A\times B\times C)&:=&\int_{[0,T]\cap
A}\,\I_{B\times C}\,[(x(t),u(t))]\,dt,
\end{eqnarray}
for all {\it rectangles} $(A\times B\times C)$, with
$(A,B,C)\in \mathcal{A}\times\mathcal{B}_{n}\times \mathcal{B}_{m}$, and 
where $\mathcal{B}_n$ (resp. $\mathcal{B}_{m}$) denotes the usual Borel $\sigma$-algebra
associated with $\RR^n$ (resp. $\RR^m$), and $\mathcal{A}$
the Borel $\sigma$-algebra of $[0,T]$, and $\I_B(\bullet)$ the indicator
function of the set $B$.

The measure $\mu ^{\mathbf{u}}$ is called the {\it occupation
measure} of the state-action (deterministic) process $(t,x(t),u(t))$ {\it up
to} time $T$, whereas $\nu ^{\mathbf{u}}$ is the occupation measure of
the state $x(T)$ {\it at} time $T$.

\begin{remark}
If the control $\mathbf{u}$ is admissible on $[0,T]$, i.e., if the
trajectory $x(\cdot)$ satisfies the constraints (\ref{cont1}) and
(\ref{cont2}), then $\nu^{\mathbf{u}}$ is a probability measure
supported on $\K$ (i.e. $\nu^\mathbf{u}\in \m(\K)_+$), and
$\mu^{\mathbf{u}}$ is supported on 
$[0,T]\times \X\times \U$ (i.e. $\mu^\mathbf{u}\in \m(\s)_+$). In particular,
$T=\mu^\mathbf{u}(\s)$.

Conversely, if the support of $\mu ^{\mathbf{u}}$ is contained in
$\s=[0,T]\times \X\times \U$ and if $\mu^\mathbf{u}(\s)=T$,
then $(x(t),u(t))\in \X\times \U$ for almost every $t\in
[0,T]$. Indeed, with (\ref{mu}),
\[T\,=\,\int_0^T
\I_{\X\times \U}\,[(x(s),u(s))]\,ds\]
\[\Rightarrow 
\I_{\X\times \U}\,[(x(s),u(s))]\,=\,1\quad \mbox{a.e. in }[0,T],\]
and hence $(x(t),u(t))\in \X\times \U$, for almost every $t\in[0,T]$.
If moreover the support of $\nu ^{\mathbf{u}}$ is
contained in $\K$, then $x(T)\in \K$. Therefore, $\mathbf{u}$ is an
admissible control on $[0,T]$.
\end{remark}

Then, observe that the optimization criterion (\ref{crit1}) now writes
$$J(0,T,x_0,\mathbf{u})\,=\,\int_{\K} H\,d\nu^{\mathbf{u}} +\int_{\s}
h\,d\mu^{\mathbf{u}}\,=\,
\langle (\mu ^{\mathbf{u}},\nu ^{\mathbf{u}}),(h,H)\rangle,$$
and one infers from (\ref{system}), (\ref{cont1}) and (\ref{cont2}) that
\begin{equation}
\label{new1}
\int_{\K} g_T \,d\nu^{\mathbf{u}} -g(0,x_0) =\int_{\s} \left(\frac{\partial g}{\partial t}+\langle \nabla_x g,f\rangle\right)
d\mu ^{\mathbf{u}},
\end{equation}
for every $g\in C_1(\Si)$ (where $g_T(x)\equiv g(T,x)$ for every $x\in
\K$), or equivalently, in view of (\ref{map}) and (\ref{map*}),
\[\langle g,\mathcal{L}^*(\mu ^{\mathbf{u}},\nu ^{\mathbf{u}})\rangle \,=\,\langle g,\delta_{(0,x_0)}\rangle,\quad \forall g\in C_1(\Si).\]
This in turn implies that
\[\mathcal{L}^*(\mu ^{\mathbf{u}},\nu ^{\mathbf{u}})\,=\, \delta_{(0,x_0)}.\]

Therefore, consider the infinite-dimensional linear program $\P$
\begin{equation}
\label{p}
\P:\quad \inf_{(\mu,\nu)\in \Delta}\:
\{\langle (\mu,\nu),(h,H)\rangle\:\vert\quad \mathcal{L}^*(\mu,\nu)=
\delta_{(0,x_0)}\}
\end{equation}
(where $\Delta:=\m(\s)_+\times \m(\K)_+$). Denote by $\inf\P$ its
optimal value and $\min\P$ is the infimum is attained, in which case
$\P$ is said to be {\it solvable}. The problem $\P$ is said
\textit{feasible} if there exists $(\mu,\nu)\in \Delta$ such that
$\mathcal{L}^*(\mu,\nu)= \delta_{(0,x_0)}$.

Note that $\P$ is feasible whenever there exists an admissible control.

The linear program $\P$ is a rephrasing
of the OCP (\ref{system})--(\ref{crit1}) in terms of the
{\it occupation measures} of its trajectories $(t,x(t),u(t))$. 
Its dual LP reads
\begin{equation}
\label{pstar}
\P^*:\quad \sup_{\varphi\in C_1(\Si)} \:\{\langle
\delta_{(0,x_0)},\varphi\rangle\:\vert
\quad \mathcal{L}\varphi \,\leq\,(h,H)\}
\end{equation}
where 
\[\mathcal{L}\varphi\leq (h,H)\Leftrightarrow\left\{\begin{array}{lll}
A\varphi(t,x,u)+h(t,x,u)&\geq 0&\forall (t,x,u)\in\s\\
\varphi(T,x)&\leq H(x)&\forall x\in \K\end{array}\right..\]
Denote by $\sup\P^*$ its optimal value and $\max\P^*$ is the supremum 
is attained, i.e. if $\P^*$ is solvable.


Discrete-time stochastic analogues of the linear programs $\P$ and $\P^*$ are also
 described in Hern\'andez-Lerma and Lasserre
\cite{ohlbook,ohlbook2} for discrete time Markov control problems. Similarly see
Cho and Stockbridge \cite{stock0}, Kurtz and Stockbridge \cite{kurtz}, and 
Helmes and Stcokbridge \cite{helmes} for some continuous-time stochastic problems.

\begin{theorem}
\label{thlp}
If $\P$ is feasible, then: 
\begin{itemize}
\item[(i)] $\P$ is solvable, i.e., $\inf\P=\min\P\leq J(0,T,x_0)$.
\item[(ii)] There is no duality gap, i.e., $\sup\P^*=\min\P$.
\item[(iii)] If moreover, for every $(t,x)\in\Si$, the set
  $f(t,x,\U)\subset\RR^n$ is convex, and the function
 \[v\mapsto g_{t,x}(v):=\inf_{u\in U}\:\{\:h(t,x,u)\::\quad v=f(t,x,u)\}\]
 is convex, then the OCP (\ref{system})--(\ref{crit1}) has an optimal
 solution and
$$\sup\P^*\,=\,\inf\P\,=\,\min\P\,=\,J^*(0,T,x_0).$$
\end{itemize}
\end{theorem}

For a proof see \S\ref{sec-thlp}. Theorem \ref{thlp}(iii) is due to Vinter \cite{vinter}.

\section{Semidefinite programming relaxations of $\P$}

The linear program $\P$ is infinite dimensional, and thus,
not tractable as it stands. Therefore, we next present a relaxation scheme
that provides a sequence of semidefinite programming, or linear matrix
inequality relaxations
(in short, LMI-relaxations) $\{\Q_r\}$,
each with {\it finitely many} constraints and variables.

Let $\RR[x]=[x_1,\ldots x_n]$ (resp. $\RR[t,x,u]=\RR[t,x_1,\ldots
x_n,u_1,\ldots,u_m]$) denote the set of polynomial functions of
the variable $x$ (resp., of the variables $t,x,u$).

Assume that $\X$ and $\K$ (resp., $\U$) are compact semi-algebraic
subsets of $\RR^n$ (resp. of $\RR^m$), of the form
\begin{eqnarray}
\label{setx}
\X&:=&\{x\in\RR^n\quad\vert\quad v_j(x)\geq 0,\quad j\in J\},\\
\label{setk}
\K&:=&\{x\in\RR^n\quad\vert\quad \theta_j(x)\geq 0,\quad j\in J_T\},\\
\label{setu}
\U&:=&\{u\in\RR^m\quad\vert\quad w_j(u)\geq 0,\quad j\in W\},
\end{eqnarray}
for some finite index sets $J_T$, $J$ and $W$, where $v_j$, $\theta_j$
and $w_j$ are polynomial functions.
Define
\begin{equation}\label{defdegXKU}
d(\X,\K,\U) := \max_{j\in J_1,\ l\in J,\ k\in W}({\rm deg}\,\theta_j,{\rm
  deg}\,v_l,{\rm deg}\,w_k).
\end{equation}

To highlight the main ideas, in this section we assume that $f$, $h$
and $H$ are polynomial functions, that is, $h\in\RR[t,x,u]$, $H\in\RR[x]$, and
$f:[0,+\infty)\times\RR^n\times \RR^m\to \RR^n$ is polynomial,
i.e., every component of $f$ satisfies $f_k\in\RR[t,x,u]$, for
$k=1,\ldots,n$.

\subsection{The underlying idea}

Observe the following important facts.

The restriction of $\RR[t,x]$ to $\Si$ belongs to $C_1(\Si)$. Therefore,
\[\mathcal{L}^*(\mu,\nu)\,=\,\delta_{(0,x_0)}\quad\Leftrightarrow\quad
\langle g,\mathcal{L}^*(\mu,\nu)\rangle \,=\,g(0,x_0),\quad\forall g\in\RR[t,x],\]
because $\Si$ being compact, polynomial functions are {\it dense} in $C_1(\Si)$
for the sup-norm. Indeed, on a compact set, one may {\it simultaneously}
approximate a function and its (continuous) partial derivatives by a
polynomial and its derivatives, uniformly (see \cite{Hirsch} pp.~65-66).
Hence, the linear program $\P$ can be written
\[\P:\left\{\begin{array}{l}
\displaystyle{\inf_{(\mu,\nu)\in \Delta}}\:\{\langle (\mu,\nu),(h,H)\rangle\\
\mbox{s.t.}\quad \langle g,\mathcal{L}^*(\mu,\nu)\rangle \,=\,g(0,x_0),\quad\forall g\in\RR[t,x],
\end{array}\right.\]
or, equivalently, by linearity,
\begin{equation}
\label{equivP}
\P:\left\{\begin{array}{l}
\displaystyle{\inf_{(\mu,\nu)\in \Delta}}\:\{\langle (\mu,\nu),(h,H)\rangle\\
\mbox{s.t.}\quad \langle \mathcal{L}g,(\mu,\nu)\rangle \,=\,g(0,x_0),\quad
\forall\, g=(t^p\,x^\alpha);\,\:(p,\alpha)\in\N\times\N^n.
\end{array}\right.
\end{equation}
The constraints of $\P$,
\begin{equation}
\label{c}
\langle \mathcal{L}g,(\mu,\nu)\rangle =g(0,x_0),\quad\forall\,
g=(t^p\,x^\alpha);\,\:(p,\alpha)\in\N\times\N^n,
\end{equation}
define countably many {\it linear} equality constraints
linking the {\it moments} of $\mu$ and
$\nu$, because if $g$ is polynomial then so 
are $\partial g/\partial t$ and $\partial g/\partial x_k$, 
for every $k$, and $\langle \nabla_x g,f\rangle$. And so, $\mathcal{L}g$
is polynomial.

The functions $h,H$ being also polynomial,
the cost $\langle (\mu,\nu),(h,H)\rangle$ of the OCP
(\ref{system})--(\ref{crit1})
is also a linear combination of the moments of $\mu$ and
$\nu$.

Therefore, the linear program $\P$ in (\ref{equivP}) can be formulated
as a LP with countably many variables (the moments of $\mu$ and $\nu$), and
countably many linear equality constraints.
However, it remains to express the fact that the variables should be moments of some measures $\mu$ and $\nu$, with support contained in $\s$ and $\K$ respectively.

At this stage, one will make some (weak) additional assumptions on the
polynomials that define the compact semi-algebraic sets $\X,\K$ and
$\U$. Under such assumptions, one may then invoke recent results of
real algebraic geometry on the representation of polynomials positive
on a compact set, and get necessary and sufficient conditions on the
variables of $\P$ to be indeed moments of two
measures $\mu$ and $\nu$, with appropriate support.
We will use Putinar's Positivstellensatz \cite{putinar} described in the next section, which yields SDP constraints on the variables.

One might also use other representation results
like e.g.\ Krivine \cite{krivine} and Vasilescu \cite{vasilescu},
and obtain {\it linear} constraints on the variables (as opposed to
SDP constraints). This is the approach taken in e.g. Helmes et al. \cite{helmes}.
However, a comparison of the use of LP-constraints versus SDP
constraints on a related problem \cite{lasserprieto} has dictated our
choice of the former.

Finally, if $g$ in (\ref{c}) runs only over all
monomials of degree less than $r$, one then obtains 
a corresponding relaxation $\Q_r$ of $\P$, which is now 
a finite-dimensional SDP that one may solve with public software
packages. At last, one lets $r\to\infty$.

\subsection{Notations, definitions and auxiliary results}
For a multi-index $\alpha=(\alpha_1,\ldots,\alpha_n)\in\N^n$, and for
$x=(x_1,\ldots,x_n)\in\RR^n$, denote $x^\alpha:=x_1^{\alpha_1}\cdots
x_n^{\alpha_n}$. Consider the canonical basis 
$\{x^\alpha\}_{\alpha\in\N^n}$ (resp., $\{t^px^\alpha
u^\beta\}_{p\in\N,\alpha\in\N^n,\beta\in\N^m}$)
of $\RR[x]$ (resp., of $\RR[t,x,u]$).

Given two sequences $y=\{y_\alpha\}_{\alpha \in\N^n}$ and
$z=\{z_{\gamma}\}_{\gamma\in\N\times\N^n\times\N^m}$ of real numbers,
define the linear functional $L_y:\RR[x]\to\RR$ by
\[H (:=\sum_{\alpha\in\N^n} H_\alpha x^\alpha )\,\mapsto L_y(H)\,:=\,
\sum_{\alpha\in\N^n} H_\alpha y_\alpha,\]
and similarly, define the linear functional $L_z:\RR[t,x,u]\to\RR$ by
\[h \mapsto L_z(h)\,:=\,\sum_{\gamma\in\N\times\N^n\times\N^m}h_\gamma\,z_\gamma\,=\,
\sum_{p\in\N,\alpha\in\N^n,\beta\in\N^m} h_{p\alpha\beta} \,z_{p\alpha\beta},\]
where $h(t,x,u)=\sum_{p\in\N,\alpha\in\N^n,\beta\in\N^m}
h_{p\alpha\beta} \,t^px^\alpha u^\beta $.

Note that, for a given measure $\nu$ (resp., $\mu$) on $\RR$ (resp.,
on $\RR\times\RR^n\times\RR^m$), there holds, for every $H\in\RR[x]$
(resp., for every $h\in \RR[t,x,u])$,
$$ \langle \nu,H\rangle = \int_\RR Hd\nu = \int_\RR \sum H_\alpha
x^\alpha d\nu = \sum H_\alpha y_\alpha = L_y(H),$$
where the real numbers $y_\alpha=\int x^\alpha d\nu$ are the moments
of the measure $\nu$ (resp., $\langle\mu,h\rangle=L_z(h)$, where $z$
is the sequence of moments of the measure $\mu$).

\begin{defn}
For a given sequence
$z=\{z_{\gamma}\}_{\gamma\in\N\times\N^n\times\N^m}$ of real numbers,
the {\it moment} matrix $M_r(z)$ of order $r$ associated with $z$, has its rows and columns indexed in the canonical basis $\{t^p x^\alpha u^\beta\}$, and is defined by \begin{equation}
\label{momentmatrix}
M_r(z)(\gamma,\beta)\,=\,z_{\gamma+\beta},\quad
\gamma,\beta\in\N\times\N^n\times\N^m,\quad
\vert\gamma\vert,\,\vert\beta\vert\leq r,
\end{equation}
where $\vert \gamma\vert :=\sum_j\gamma_j$.
The moment matrix $M_r(y)$ of order $r$ associated with a given
sequence $y=\{y_\alpha\}_{\alpha \in\N^n}$, has its rows and columns indexed in the canonical basis $\{x^\alpha\}$, and is defined in a similar fashion.
\end{defn}

Note that, if $z$ has a representing measure $\mu$, i.e., if $z$ is the
sequence of moments of the measure $\mu$ on
$\RR\times\RR^n\times\RR^m$, then $L_z(h)=\int hd\mu$, for every
$h\in\RR[t,x,u]$, and if $\mathbf{h}$ denotes the vector of
coefficients of $h\in\RR[t,x,u]$ of degree less than $r$, then
\[\langle \mathbf{h},M_r(z)\mathbf{h}\rangle\,=\,L_z(h^2)\,=\,\int
h^2\,d\mu\geq0.\]
This implies that $M_r(z)$ is symmetric nonnegative (denoted
$M_r(z)\succeq0$), for every $r$. The same holds for $M_r(y)$.

Conversely, not every sequence $y$ that satisfies $M_r(y)\succeq0$ for
every $r$, has a representing measure. However, several sufficient conditions
exist, and in particular the following one, due to Berg \cite{Berg}.

\begin{proposition}
\label{prop1}
If $y=\{y_\alpha\}_{\alpha\in\N^n}$ satisfies $\vert y_\alpha\vert\leq
1$ for every $\alpha\in\N^n$, and $M_r(y)\succeq0$ for every
integer $r$, then $y$ has a representing measure on
$\RR^n$, with support contained in the unit ball $[-1,1]^{n}$.
\end{proposition}

We next present another sufficient condition which is crucial in the
proof of our main result.

\begin{defn}
For a given polynomial $\theta\in\RR[t,x,u]$, written
$$\theta(t,x,u)=\sum_{\delta=(p,\alpha,\beta)}\theta_\delta\, t^px^\alpha
u^\beta,$$
define the {\it localizing} matrix $M_r(\theta\,z)$
associated with $z,\theta$, and with rows and columns also indexed in
the canonical basis of $\RR[t,x,u]$, by
\begin{equation}
\label{localizingmatrix}
M_r(\theta\,z)(\gamma,\beta)\,=\,\sum_{\delta}\theta_\delta \,z_{\delta+\gamma+\beta}
\quad \gamma,\beta\in\N\times\N^n\times\N^m,\quad \vert\gamma\vert,\,\vert\beta\vert\leq r.
\end{equation}
The localizing matrix $M_r(\theta\,y)$ associated with a given
sequence $y=\{y_\alpha\}_{\alpha \in\N^n}$ is defined similarly.
\end{defn}

Note that, if
$z$ has a representing measure $\mu$ on $\RR\times\RR^n\times\RR^m$
with support contained in the level set $\{(t,x,u) :\:
\theta(t,x,u)\geq 0\}$, and if $h\in\RR[t,x,u]$ has degree less than
$r$, then
\[\langle
\mathbf{h},M_r(\theta,z)\mathbf{h}\rangle\,=\,L_z(\theta\,h^2)\,=\,\int
\theta h^2\,d\mu\geq 0.\]
Hence, $M_r(\theta\,z)\succeq 0$, for every $r$.

Let $\Sigma^2\subset\RR[x]$ be the convex cone generated in $\RR[x]$
by all squares of polynomials,
and let $\Omega\subset\RR^n$ be the compact basic semi-algebraic set defined by
\begin{equation}
\label{setomega}
\Omega\,:=\,\{x\in\RR^n\quad\vert\quad g_j(x)\geq0,\quad j=1,\ldots,m\}
\end{equation}
for some family $\{g_j\}_{j=1}^m\subset\RR[x]$.

\begin{defn}
\label{def1}
The set $\Omega\subset\RR^n$ defined by 
(\ref{setomega}) satisfies Putinar's condition if
there exists $u\in\RR[x]$ such that $u=u_0+\sum_{j=1}^mu_jg_j$ for some family 
$\{u_j\}_{j=0}^m\subset\Sigma^2$, and the level set
$\{x\in\RR^n\ \vert\ u(x)\geq0\}$ is compact.
\end{defn}
Putinar's condition is satisfied if e.g. the level set 
$\{x: \:g_k(x)\geq0\}$ is compact for some $k$, or if all the $g_j$'s are linear, in which case $\Omega$ is a polytope. In addition, if one knows some $M$ such that $\Vert x\Vert\leq M$ 
whenever $x\in\Omega$, then it suffices to add the redundant quadratic constraint 
$M^2-\Vert x\Vert^2\geq0$ in the definition (\ref{setomega}) of $\Omega$, and Putinar's condition is satisfied (take $u:=M^2-\Vert x\Vert^2$).

\begin{theorem}[Putinar's Positivstellensatz \cite{putinar}]
\label{putinarthm}
Assume that the set $\Omega$ defined by (\ref{setomega})
satisfies Putinar's condition.
\begin{itemize}
\item[(a)] If $f\in\RR[x]$ and $f>0$ on $\Omega$, then
\begin{equation}
\label{puta}
f=f_0+\sum_{j=1}^mf_j\,g_j,
\end{equation}
for some family $\{f_j\}_{j=0}^m\subset\Sigma^2$.
\item[(b)] Let $y=\{y_\alpha\}_{\alpha\in\N^n}$ be a sequence of real
  numbers. If
\begin{equation}
\label{putb}
M_r(y)\succeq0\,;\quad M_r(g_j\,y)\succeq0,\quad j=1,\ldots,m;\quad\forall\,r=0,1,\ldots
\end{equation}
then $y$ has a representing measure with support contained in $\Omega$.
\end{itemize}
\end{theorem}

\subsection{LMI-relaxations of $\P$}
Consider the linear program $\P$ defined by (\ref{equivP}).

Since the semi-algebraic sets $\X,\K$ and $\U$ defined respectively by
(\ref{setx}), (\ref{setk}) and (\ref{setu}) are compact, with no loss
of generality, we assume (up to a scaling of the variables $x,u$ and $t$) that 
$T=1$, $\X,\K\subseteq [-1,1]^n$ and $\U\subseteq [-1,1]^m$. 

Next, given a sequence $z=\{z_\gamma\}$ indexed in the basis of $\RR[t,x,u]$ denote
$z(t)$, $z(x)$ and $z(u)$ its marginals with respect to the variables $t$, $x$ and
$u$, respectively. These sequences are indexed in the canonical basis of
$\RR[t]$, $\RR[x]$ and $\RR[u]$ repectively.
For instance, writing
$\gamma=(k,\alpha,\beta)\in\N\times\N^n\times\N^n$,
\[\{z(t)\}=\{z_{k,0,0}\}_{k\in\N};\quad
\{z(x)\}=\{z_{0,\alpha,0}\}_{\alpha\in\N^n};\quad
\{z(u)\}=\{z_{0,0,\beta}\}_{\beta\in\N^m}.\]

Let $r_0$ be an integer such that
$2r_0\geq \max\,({\rm deg}\,f,{\rm deg}\,h,{\rm
  deg}\,H,\,2d(\X,\K,\U) )$, where $d(\X,\K,\U)$ is defined by
(\ref{defdegXKU}).
For every $r\geq r_0$, consider
the LMI-relaxation
\begin{equation}
\label{sdpq}
\Q_r:\left\{\begin{array}{l}
\displaystyle{\inf_{y,z}}\quad L_z(h)+L_y(H)\\
M_r(y),\,M_r(z)\succeq 0\\
M_{r-{\rm deg}\,\theta_j}(\theta_j\,y)\succeq0,\quad j\in J_1\\
M_{r-{\rm deg}\,v_j}(v_j\,z(x))\succeq0,\quad j\in J\\
M_{r-{\rm deg}\,w_k}(w_k\,z(u))\succeq0,\quad k\in W\\
M_{r-1}(t(1-t)\,z(t))\succeq 0\\
L_y(g_1) -L_z(\partial g/\partial t+\langle \nabla_x
g,f\rangle)=g(0,x_0),\quad \forall g=(t^px^\alpha)\\
\mbox{with }\:p+\vert\alpha\vert-1+{\rm deg}\,f\leq 2r\\
\end{array}\right.,
\end{equation}
whose optimal value is denoted by $\inf\Q_r$. 

\vspace{0.2cm}

\noindent
{\bf OCP with free terminal time.}
For the OCP (\ref{crit1libre}), i.e., with {\it free} terminal time
$T\leq T_0$, we need adapt the notation because $T$ is now a variable.
As already mentioned in Remark \ref{continuity}(iii), the measure
$\nu$ in the infinite dimensional linear program $\P$ defined in (\ref{p}), is now supported in 
$[0,T_0]\times \K$ (and $[0,1]\times \K$ after re-scaling) instead of $\K$ previously.
Hence, the sequence $y$ associated with $\nu$ is now indexed in the basis 
$\{t^px^\alpha\}$ of $\RR[t,x]$ instead of $\{x^\alpha\}$ previously. Therefore, given
$y=\{y_{k\alpha}\}$ indexed in that basis, let $y(t)$ and $y(x)$ be the subsequences of $y$ defined by:
\[y(t)\,:=\,\{y_{k0}\}_{k},\quad k\in\N;;\quad y(x)\,=\,\{y_{0\alpha}\},\quad \alpha\in\N^n.\]

Then again (after rescaling), the LMI-relaxation $\Q_r$ now reads
\begin{equation}
\label{sdpqrt}
\Q_r:\left\{\begin{array}{l}
\displaystyle{\inf_{y,z}}\quad L_z(h)+L_y(H)\\
M_r(y),\,M_r(z)\succeq0\\
M_{r-r(\theta_j)}(\theta_j\,y)\succeq0,\quad j\in J_1\\
M_{r-r(v_j)}(v_j\,z(x))\succeq0,\quad j\in J\\
M_{r-r(w_k)}(w_k\,z(u))\succeq0,\quad k\in W\\
M_{r-1}(t(1-t)\,y(t))\succeq 0\\
M_{r-1}(t(1-t)\,z(t))\succeq 0\\
L_y(g) -L_z(\partial g/\partial t+\langle \nabla_x g,f\rangle)=g(0,x_0),\quad
\forall g=(t^px^\alpha)\\
\mbox{with }\:p+\vert\alpha\vert-1+{\rm deg}\,f\leq 2r\\
\end{array}\right..
\end{equation}
The particular case of minimal time problem is obtained with
$h\equiv 1$, $H\equiv 0$.

For {\it time-homogeneous} problems, i.e., when $h$ and $f$ do not depend on $t$,
one may take $\mu$ (resp. $\nu$) supported on $\X\times\U$
(resp. $\K$), which simplifies the associated LMI-relaxation (\ref{sdpqrt}).

The main result is the following.

\begin{theorem}
\label{main}
Let $\X,\K\subset [-1,1]^n$, and $\U\subset [-1,1]^m$
be compact 
basic semi-algebraic-sets respectively defined by (\ref{setx}), (\ref{setk})
and (\ref{setu}).
Assume that $\X,\K$ and $\U$ satisfy Putinar's condition (see
Definition (\ref{def1})), and
let $\Q_r$ be the LMI-relaxation defined in (\ref{sdpq}). Then,
\begin{itemize}
\item[(i)] $\inf\Q_r\uparrow \min\P$ as $r\to\infty$;
\item[(ii)] if moreover, for every $(t,x)\in\Si$, the set
  $f(t,x,\U)\subset\RR^n$ is convex, and the function
 \[v\mapsto g_{t,x}(v):=\inf_{u\in U}\:\{\:h(t,x,u)\ \vert\ v=f(t,x,u)\}\]
 is convex, then $\inf\Q_r\uparrow \min\P=J^*(0,T,x_0)$, as
 $r\to\infty$.
\end{itemize}
\end{theorem}

The proof of this result is postponed to the Appendix in Section \S
\ref{proofmain}.

\subsection{The dual $\Q^*_r$}

We describe the dual of the LMI-relaxation $\Q_r$ which is also 
a  semidefinite program, denoted
$\Q^*_r$, and relate $\Q^*_r$ with the dual $\P^*$ of $\P$, defined in
(\ref{pstar}).\\

Let $s(r)$ be the cardinal of the set
$V_r:=\{(k,\alpha)\in\N\times\N^n\ \vert\ k+\vert\alpha\vert\leq
r-r_0\}$, and given $\lambda\in\RR^{s(r)}$, let $\Lambda_{r}\in\RR[t,x]$ be the polynomial
\[(t,x)\,\mapsto\,\Lambda_{r}(t,x)\,:=\,\sum_{(k,\alpha)\in
  V_r}\lambda_{k\alpha}\,t^kx^\alpha.\] 
Consider the semidefinite program:
\begin{equation}
\label{sdpqdual}
\Q^*_r:\left\{\begin{array}{l}
\displaystyle{\sup_{q_0,q_j^x,q_k^y,l_0,l_j,\Lambda_r}}\quad
\Lambda_{ r}(0,x_0),\\ 
\\
h+A\,\Lambda_{r}\,=\,q_0\,t(1-t)+\sum_{k\in W}q_k^u\,w_k
+\sum_{j\in J}q_j^x\,v_j,\\
\\
H-\Lambda_{r}(1,.)\,=\,l_0+\sum_{j\in J_1}l_j\,\theta_j,\\
\\
q_0\in\RR[t],\,q^u_k\in\RR[u],\,q^x_j\in\RR[x],\,l_j\in\RR[x]\\
\\
\{q_0,q_j^x,q_k^u,l_0,l_j\}\mbox{ s.o.s. (sums of squares
polynomials)}, \quad\mbox{and}\\ 
{\rm deg}\,l_j\theta_j,\:
{\rm deg}\,q_j^xv_j,\:
{\rm deg}\,q_k^uw_k,\:
{\rm deg}\,q_0\leq\,2r-2;\:{\rm deg}\,l_0\,\leq\,2r.
 \end{array}\right..
\end{equation}
The LMI $\Q^*_r$ is a reinforcement of $\P^*$ in the following sense:
\begin{itemize}
\item the unknown
function $\varphi\in C_1(\Si)$ is now replaced with a polynomial
$\Lambda_{r}\in\RR[t,x]$ of degree less than $2r$;
\item the constraint $-A\varphi \leq h$ for $(t,x,u)\in\s$, is now
  replaced with the constraint $h+A\Lambda_{ r}\geq0$ on $\s$
  and the polynomial $h+A\Lambda_{r}\geq0$ which is nonnegative
  on $\s$, has Putinar's representation $q_0\,t(1-t)+\sum_{k\in
    W}q_k^u\,w_k +\sum_{j\in J}q_j^x\,v_j$;
\item the constraint $\varphi_1\leq H$ for $x\in\K$, is replaced
  with the constraint $H-\Lambda_{r}(1,.)\geq0$ on $\K$, and
  the polynomial $H-\Lambda_{r}(1,.)$ which is nonnegative on
  $\K$, has Putinar's representation $l_0+\sum_{j\in
    J_1}l_j\,\theta_j$.
\end{itemize}

Assume that $\Q^*_r$ is solvable. A natural question is to know
whether or not we can use an optimal solution
\[q_0,q_j^x,q_k^y,l_0,l_j,\Lambda_r\] of $\Q^*_r$ to obtain some
information on an optimal solution of $\P$.
The most natural idea is to look for the zero set in $\s$ of the polynomial
\[(t,x,u)\mapsto\,q_0\,t(1-t)+\sum_{k\in W}q_k^uw_k
+\sum_{j\in J}q_j^xv_j.\]
Indeed, under the assumptions of Theorem \ref{main}, if
$\sup\Q^*_r=\inf\Q_r$, then
$\Lambda_{r}(0,x_0)\approx\inf\Q_r\approx\min\P=\sup\P^*$, and so, the polynomial $\Lambda_{r}\in\RR[t,x]$ seems to be a good candidate 
to approximate a nearly optimal solution $\varphi\in C_1(\Si)$ of $\P^*$.

Next, as $\Q^*_r$ is an approximation of a
weak formulation of the HJB optimality equation,
one may hope that the zero set in $\s$ 
of the polynomial $h+A\Lambda_{ r}$ provides some good information
on the possible states $x^*(t)$ and controls $u^*(t)$ at time $t$
in an optimal solution of the OCP (\ref{system})--(\ref{crit1}).

That is, fixing an arbitrary $t_0\in[0,1]$, one may solve the equation
\[\sum_{k\in W}q_k^u(u)\,w_k(u)+\sum_{j\in J}q_j^x(x)\,v_j(x)\,=\,-q_0\,t_0(1-t_0),\]
and look for solutions $(x,u)\in \X\times \U$.

All these issues deserve further investigation beyond the scope of the present paper.  However, at least in the minimum time 
problem for the (state and control constrained) double integrator example considered in \S\ref{ex1},
we already have some numerical support for the above claims.

\subsection{Certificates of non controllability}

For minimum time OCPs, i.e., with free terminal time $T$ and instantaneous cost 
$h\equiv 1$, and $H\equiv 0$, the LMI-relaxations $\Q_r$ defined in (\ref{sdpqrt})
may provide {\it certificates} of non controllability. 

Indeed, if for a given initial state $x_0\in\X$, some LMI-relaxation $\Q_r$ in the hierarchy has {\it no} feasible solution, then the initial state $x_0$ {\it cannot} be steered to the origin in finite time.
In other words, $\inf\Q_r=+\infty$ provides a certificate of uncontrollability 
of the initial state $x_0$.  
It turns out that sometimes such certificates can be provided at cheap cost, i.e., 
with LMI-relaxations of low order $r$. This is illustrated on the Zermelo problem
in \S\ref{controllability}.

Moreover, one may also consider controllability in given finite time $T$, that is, 
consider the LMI-relaxations as defined in (\ref{sdpq}) with $T$ fixed, and 
$H\equiv 0$, $h\equiv 1$. Again,
if for a given initial state $x_0\in\X$, the LMI-relaxation $\Q_r$ has no feasible solution, the initial state $x_0$ cannot be steered to the origin in less than $T$ units of time.
And so, $\inf\Q_r=+\infty$ also provides a certificate of uncontrollability 
of the initial state $x_0$.  

\section{Generalization to smooth optimal control problems}
\label{smooth}

In the previous section, we assumed, to highlight the main ideas, that
$f$, $h$ and $H$ were polynomials. In this section, we generalize
Theorem \ref{main}, and simply assume that $f$, $h$ and $H$ are smooth.
Consider the linear program $\P$ defined in the previous section
$$\P:\left\{\begin{array}{l}
\displaystyle{\inf_{(\mu,\nu)\in \Delta}}\:\{\langle
(\mu,\nu),(h,H)\rangle\\
\mbox{s.t.}\quad \langle g,\mathcal{L}^*(\mu,\nu)\rangle
\,=\,g(0,x_0),\quad\forall g\in\RR[t,x].
\end{array}\right.$$

Since the sets $\X$, $\K$ and $\U$, defined
previously, are compact, it follows from \cite{Coatmelec} (see also
\cite[pp.~65-66]{Hirsch}) that $f$ (resp.\ $h$, resp.\ $H$) is the limit
in $C_1(\s)$ (resp.\ $C_1(\s)$, resp.\ $C_1(\K)$) of a sequence of
polynomials $f_p$ (resp.\ $h_p$, resp.\ $H_p$) of degree $p$, as
$p\rightarrow +\infty$.

Hence, for every integer $p$, consider
the linear program $\P_p$
$$\P_p:\left\{\begin{array}{l}
\displaystyle{\inf_{(\mu,\nu)\in \Delta}}\:\{\langle
(\mu,\nu),(h_p,H_p)\rangle\\
\mbox{s.t.}\quad \langle g,\mathcal{L}_p^*(\mu,\nu)\rangle
\,=\,g(0,x_0),\quad\forall g\in\RR[t,x],
\end{array}\right.$$
the smooth analogue of $\P$,
where the linear mapping $\mathcal{L}_p:\,C_1(\Si)\to C(\s)\times
C(\K)$ is defined by
$$\mathcal{L}_p\varphi \,:=\,(-A_p\varphi, \varphi_T),$$
and where $A_p:\,C_1(\Si)\to C(\s)$ is defined by
$$
A_p\varphi (t,x,u)\,:=\,\frac{\partial \varphi}{\partial
t}(t,x)+\langle f_p(t,x,u),\nabla_x \varphi(t,x)\rangle.
$$
For every integer $r\geq \max( p/2 , d(\X,\K,\U))$,
let $\Q_{r,p}$ denote the LMI-relaxation
(\ref{sdpq}) associated with the linear program $\P_p$.

Recall that from Theorem \ref{main}, if $\K$, $\X$ and $\U$ satisfy Putinar's
condition, then $\inf\Q_{r,p}\uparrow \min\P_p$ as $r\to +\infty$;

The next result, generalizing Theorem \ref{main}, shows that
it is possible to let $p$ tend to $+\infty$.
For convenience, set
$$v_{r,p}=\inf\Q_{r,p},\ \ v_p=\min\P_p,\ \ v=\min\P.$$

\begin{theorem}
\label{mainbis}
Let $\X,\K\subset [-1,1]^n$, and $\U\subset [-1,1]^m$
be compact semi-algebraic-sets respectively defined by (\ref{setx}),
(\ref{setk}) and (\ref{setu}).
Assume that $\X,\K$ and $\U$ satisfy Putinar's condition (see
Definition (\ref{def1})). Then,
\begin{itemize}
\item[(i)] $\displaystyle
  v=\lim_{p\to+\infty}\lim_{\underset{\scriptstyle
      2r>p}{r\to+\infty}}v_{r,p}
= \lim_{p\to+\infty}\sup_{r>p/2}v_{r,p}\ \leq\ J^*(0,T,x_0)$.
\item[(ii)] Moreover if for every $(t,x)\in\Si$, 
the set $f(t,x,\U)\subset\RR^n$ is convex,
and the function
 \[v\mapsto g_{t,x}(v):=\inf_{u\in U}\:\{\:h(t,x,u)\ \vert\ v=f(t,x,u)\}\]
is convex, then $v=J^*(0,T,x_0)$.
\end{itemize}
\end{theorem}

The proof of this result is in the Appendix, Section \S
\ref{proofmainbis}.

From the numerical point of
view, depending on the functions $f$, $h$, $H$, the degree of the polynomials
of the approximate OCP $\P_p$ may be required to be large, and hence the
hierarchy of LMI-relaxations  $(\Q_r)$ in (\ref{sdpq}) 
might not be efficiently implementable, at least
in view of the performances of public SDP solvers available at
present.


\begin{remark}
The previous construction extends to smooth optimal control problems
on Riemannian manifolds, as follows.
Let $M$ and $N$ be smooth Riemannian manifolds. Consider on $M$ the
control system (\ref{system}),
where $f:[0,+\infty)\times M \times N \longrightarrow TM$ is smooth,
and where the controls are bounded measurable functions, defined on intervals
$[0,T(\mathbf{u})]$ of $\RR^+$, and taking their values in a
\textit{compact} subset $\U$ of $N$.
Let $x_0\in M$, and let $\X$ and $\K$ be compact subsets of $M$.
Admissible controls are defined as previously.
For an admissible control $\mathbf{u}$ on $[0,T]$, the cost of the
associated trajectory $x(\cdot)$ is defined by (\ref{cost}), where
where $h:[0,+\infty)\times M \times N \longrightarrow \RR$
and $H:M\rightarrow \RR$ are smooth functions.

According to Nash embedding Theorem \cite{Nash}, there exist an
integer $n$ (resp.\ $m$) such that $M$ (resp.\ $N$) is smoothly
isometrically embedded in $\RR^n$ (resp.\ $\RR^m$). In this context,
all previous results apply.

This remark is important for the applicability of the method described
in this article. Indeed, many practical control problems (in
particular, in mechanics) are expressed
on manifolds, and since the optimal control problem investigated here
is global, they cannot be expressed in general as control systems in
$\RR^n$ (in a global chart).
\end{remark}

\section{Illustrative examples}

We consider here the minimal time OCP, that is, we aim
to approximate the {\it minimal time} to steer 
a given initial condition to the origin. We have tested  
the above methodology on two test OCPs,
the double and Brockett {\it integrators}, for which the associated
optimal value $T ^*$ can be calculated exactly.
The numerical examples in this section were processed with
our Matlab package \textit{GloptiPoly 3}\footnote{
\textit{GloptiPoly 3} can be downloaded at
{\tt www.laas.fr/$\sim$henrion/software}}.

\subsection{The double integrator}
\label{ex1}
Consider the double integrator system in $\RR^2$
\begin{equation}\label{inte}
\begin{array}{rcl}
\dot{x}_1(t)&=&x_2(t),
\\
\dot{x}_2(t)&=&u(t),
\end{array}
\end{equation}
where $x=(x_1,x_2)$ is the state and the control ${\mathbf u}=u(t)\in\mathcal{U}$,
satisfies the constraint $\vert u(t)\vert\leq 1$, for all $t\geq 0$. In addition, the state is constrained to satisfy $x_2(t)\geq-1$, for all $t$.
In this time-homogeneous case, and with the notation of Section
\ref{general},  we have $\X=\{x\in\RR^2\::\:x_2\geq-1\}$, $\K =\{(0,0)\}$, and $\U=[-1,1]$.

\begin{remark}
The theorem obviously extends, up to scaling, to the case of
arbitrary compact subsets $\X,\K\subset\RR^n$ and $\U\subset[-1,1]^m$.
\end{remark}

Observe that $\X$ is not compact and so the convergence result of Theorem \ref{main} may not hold. In fact, we may impose the additional constraint $\Vert x(t)\Vert_\infty\leq M$ for some large $M$ 
(and modify $\X$ accordingly), because for initial states $x_0$ with 
$\Vert x_0\Vert_\infty$ relatively small with respect to $M$, the optimal trajectory remains in $\X$.
However, in the numerical experiments, we have not enforced an additional constraint.
We have maintained the original constraint $x_2\geq-1$ in the localizing constraint $M_{r-r(v_j)}(v_jz(x))\succeq0$, with $x\mapsto v_j(x)=x_2+1$.

\subsubsection{Exact computation}

For this very simple system, one is able to compute exactly the
optimal minimum time.
Denoting
$T(x)$ the minimal time to reach the origin from
$x=(x_1,x_2)$, we have:

If $x_1 \geq 1-x_2^2/2$ and $x_2 \geq -1$ then
$T(x)=x_2^2/2+x_1+x_2+1$.
If $-x_2^2/2\,{\rm sign}\,x_2 \leq x_1 \leq 1-x_2^2/2$
and $x_2 \geq -1$ then 
$T(x)= 2\sqrt{x_2^2/2+x_1}+x_2$.
If $x_1 < -x_2^2/2\,{\rm sign}\,x_2$ and $x_2 \geq -1$
then $T(x)= 2\sqrt{x_2^2/2-x_1}-x_2$.

\subsubsection{Numerical approximation}

Table \ref{doubleinte-x0} displays the values of the initial state $x_0\in\X$, and 
denoting $\inf\Q_r(x_0)$ the optimal value of the LMI-relaxation (\ref{sdpqrt}) for
the minimum time OCP (\ref{inte}) with initial state $x_0$,
Tables \ref{doubleinte-2},  \ref{doubleinte-3}, and  \ref{doubleinte-5}
display the numerical values of the ratii $\inf\Q_r(x_0)/T(x_0)$ for $r=2,3$ and $5$ respectively.

{\footnotesize\begin{table}
\caption{Double integrator: data initial state $x_0=(x_{01},x_{02})$}
\label{doubleinte-x0}
\begin{center}
\begin{tabular}{||c|c|c|c|c|c|c|c|c|c|c|c||}
\hline
\hline
$x_{01}$  & 0.0   & 0.2 &   0.4 &   0.6 &   0.8   & 1.0  &  1.2&    1.4&    1.6  &  1.8  &  2.0\\
\hline
 $x_{02}$ & -1.0 &  -0.8 &  -0.6&   -0.4 &  -0.2    &     0.0   & 0.2  &  0.4  &  0.6  &  0.8  &  1.0 \\ 
\hline
    \hline
\end{tabular}
\end{center}
\end{table}}

{\footnotesize\begin{table}
\caption{Double integrator: ratio $\inf\Q_2/T(x_0)$}
\label{doubleinte-2}
\begin{center}
\begin{tabular}{||c|c|c|c|c|c|c|c|c|c|c||}
\hline
\multicolumn{11}{||c||}{second LMI-relaxation: r=2}\\
\hline
\hline
\hline
  0.4598 &   0.3964 &   0.3512   & 0.9817&   0.9674 &   0.9634  &  0.9628  &  0.9608  &  0.9600 &   0.9596  &  0.9595\\
\hline
    0.4534  &  0.3679&    0.9653 &   0.9347 &   0.9355 &   0.9383  &  0.9385  &  0.9386 &   0.9413  &  0.9432  &  0.9445\\
\hline
    0.4390  &  0.9722  &  0.8653  &  0.8457  &  0.8518  &  0.8639  &  0.8720  &  0.8848  &  0.8862  &  0.8983   & 0.9015\\
\hline
    0.4301  &  0.7698   & 0.7057 &   0.7050 &   0.7286  &  0.7542  &  0.7752   & 0.7964  &  0.8085   & 0.8187 &   0.8351\\
\hline
    0.4212   & 0.4919  &  0.5039  &  0.5422  &  0.5833  &  0.6230  &  0.6613  &  0.6870 &   0.7121   & 0.7329   & 0.7513\\
\hline
    0.0000  &  0.2238 &   0.3165 &  0.3877 &   0.4476 &   0.5005   & 0.5460   & 0.5839  &  0.6158 &   0.6434  &  0.6671\\
\hline
    0.4501  &  0.3536 &   0.3950  &  0.4403  &  0.4846  &  0.5276 &   0.5663   & 0.5934 &   0.6204 &   0.6474   & 0.6667\\
\hline
    0.4878   & 0.4493 &   0.4699  &  0.5025  &  0.5342 &   0.5691  &  0.5981  &  0.6219  &  0.6446 &   0.6647  &  0.6824\\
\hline
    0.5248   & 0.5142  &  0.5355   & 0.5591  &  0.5840  &  0.6124  &  0.6312   & 0.6544  &  0.6689 &   0.6869  &  0.7005\\
\hline
    0.5629  &  0.5673   & 0.5842  &  0.6044  &  0.6296 &   0.6465 &   0.6674 &   0.6829  &  0.6906 &   0.7083  &  0.7204\\
\hline
    0.6001  &  0.6099  &  0.6245  &  0.6470 &   0.6617  &  0.6792 &   0.6972   & 0.7028  & 0.7153   & 0.7279  &  0.7369\\
    \hline
    \hline
\end{tabular}
\end{center}
\end{table}}

{\footnotesize\begin{table}
\caption{Double integrator: ratio $\inf\Q_3/T(x_0)$}
\label{doubleinte-3}
\begin{center}
\begin{tabular}{||c|c|c|c|c|c|c|c|c|c|c||}
\hline
\multicolumn{11}{||c||}{third LMI-relaxation: r=3}\\
\hline
\hline
 0.5418  &  0.4400   & 0.3630  &  0.9989  &  0.9987  &  0.9987  &  0.9985  &  0.9984  &  0.9983  &  0.9984  &  0.9984\\
 \hline
    0.5115   & 0.3864  &  0.9803 &   0.9648 &   0.9687  &  0.9726  &  0.9756   & 0.9778 &   0.9798 &   0.9815  &  0.9829\\
    \hline
    0.4848  &  0.9793 &  0.8877  &  0.8745  &  0.8847   & 0.8997  &  0.9110  &  0.9208  &  0.9277   & 0.9339 &   0.9385\\
 \hline
     0.4613  &  0.7899  &  0.7321  &  0.7401 &   0.7636   & 0.7915  &  0.8147  &  0.8339   & 0.8484  &  0.8605  &  0.8714\\
  \hline
     0.4359   & 0.5179   & 0.5361 &   0.5772  &  0.6205  &  0.6629  &  0.7013 &   0.7302&    0.7540   & 0.7711   & 0.7891\\
   \hline
     0.0000  &  0.2458   & 0.3496&    0.4273&    0.4979   & 0.5571  &  0.5978  &  0.6409 &   0.6719  &  0.6925  &  0.7254\\
    \hline
     0.4556  &  0.3740  &  0.4242  &  0.4789   & 0.5253  &  0.5767 &   0.6166  &  0.6437  &  0.6807 &   0.6972   & 0.7342\\
     \hline
     0.4978  &  0.4709   & 0.5020  &  0.5393  &  0.5784 &   0.6179   & 0.6477 &   0.6776  &  0.6976  &  0.7192  &  0.7347\\
     \hline
     0.5396 &   0.5395 &   0.5638 &  0.5955   & 0.6314   & 0.6600  &  0.6856  &  0.7089  &  0.7269  &  0.7438  &  0.7555\\
     \hline
     0.5823  &  0.5946  &  0.6190  &  0.6453   & 0.6703  &  0.7019  &  0.7177  &  0.7382  &  0.7539  &  0.7662  &  0.7767\\
     \hline
     0.6255  &  0.6434  &  0.6656  &  0.6903  &  0.7193  &  0.7326   & 0.7543   & 0.7649  &  0.7776  &  0.7917 &   0.8012\\
     \hline
    \hline
\end{tabular}
\end{center}
\end{table}}

\begin{figure}
\resizebox{0.67\textwidth}{!}
{\includegraphics{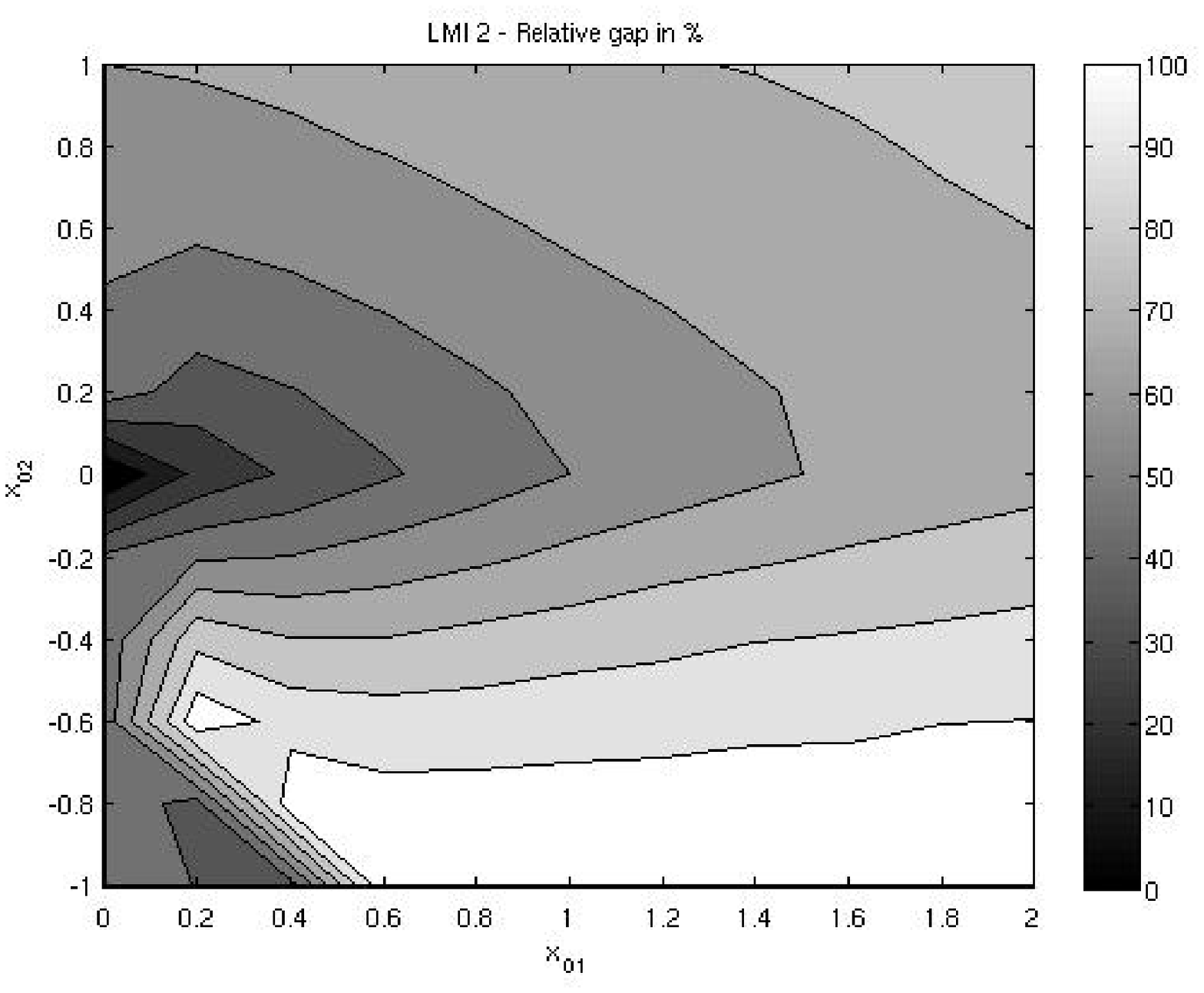}}
\caption{Double integrator: level sets $\inf\Q_2/T(x_0)$}
\label{figinte2}
\end{figure}

\begin{figure}
\resizebox{0.67\textwidth}{!}
{\includegraphics{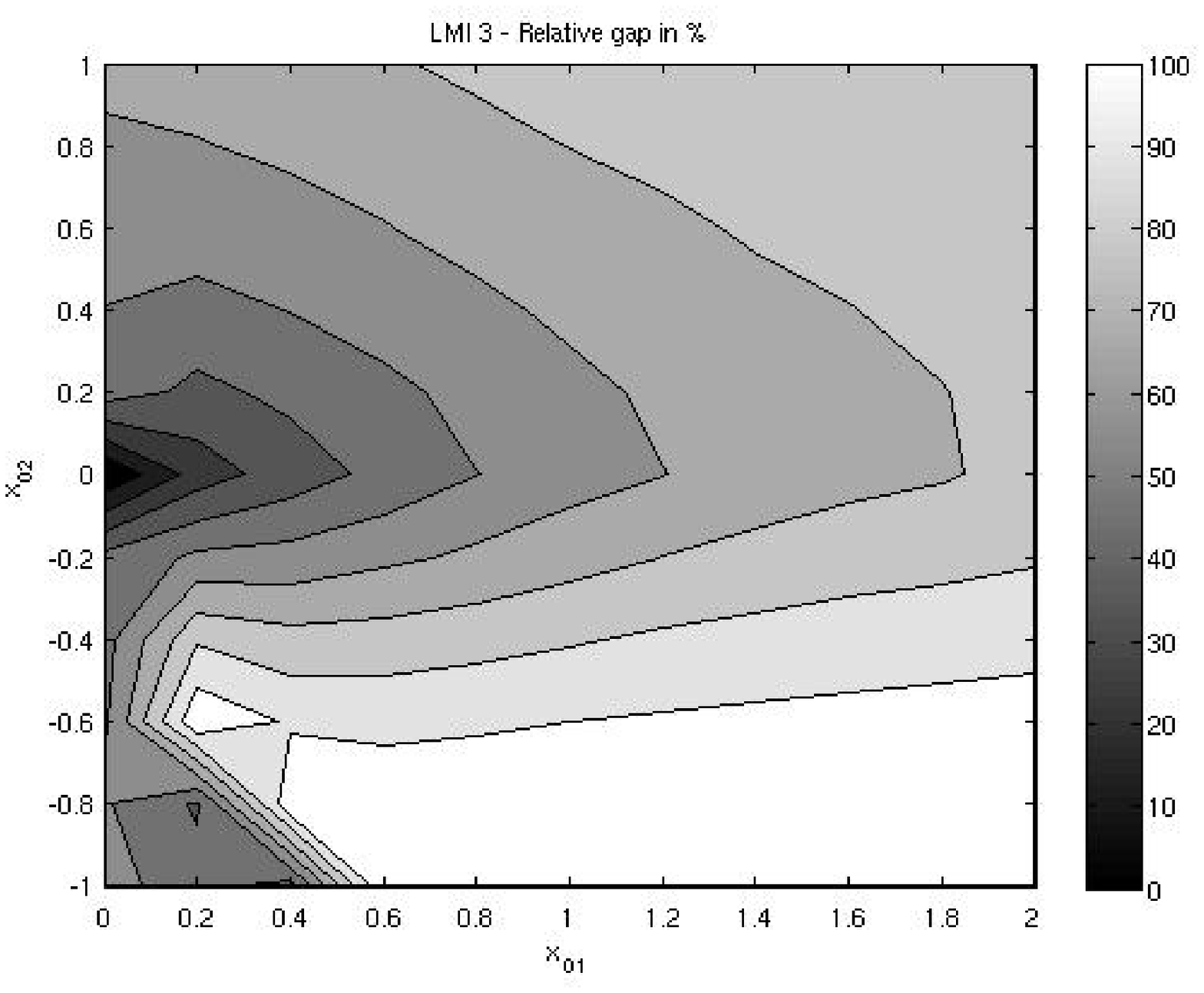}}
\caption{Double integrator: level sets $\inf\Q_3/T(x_0)$}
\label{figinte3}
\end{figure}

\begin{figure}
\resizebox{0.67\textwidth}{!}
{\includegraphics{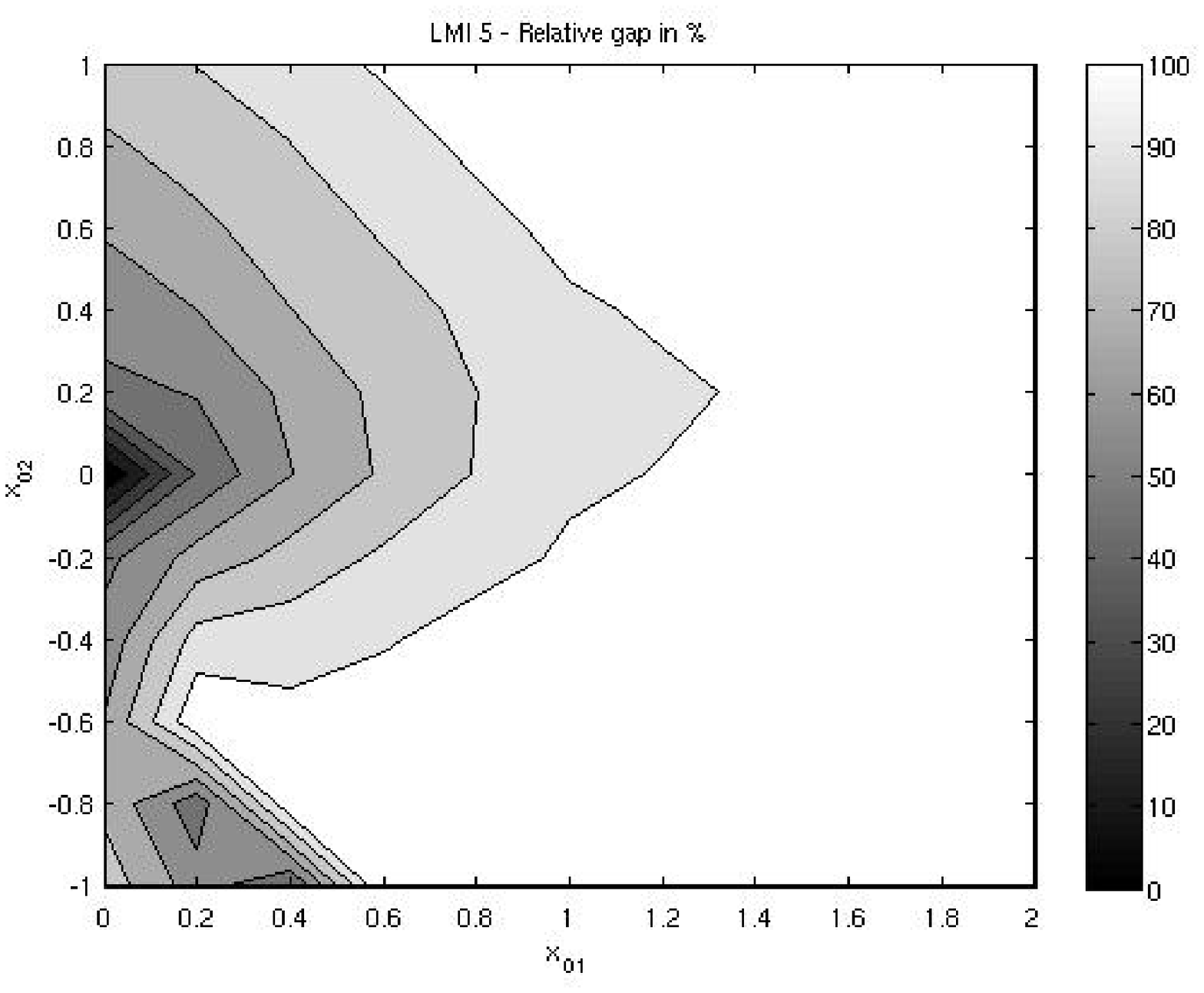}}
\caption{Double integrator: level sets $\inf\Q_5/T(x_0)$}
\label{figinte5}
\end{figure}

{\footnotesize\begin{table}
\caption{Double integrator: ratio $\inf\Q_5/T(x_0)$}
\label{doubleinte-5}
\begin{center}
\begin{tabular}{||c|c|c|c|c|c|c|c|c|c|c||}
\hline
\multicolumn{11}{||c||}{fifth LMI-relaxation: r=5}\\
\hline
\hline
 0.7550  &  0.5539  &  0.3928 &   0.9995  &  0.9995  &  0.9995  &  0.9994  &  0.9992  &  0.9988  &  0.9985   & 0.9984\\
 \hline
0.6799 &   0.4354 &   0.9828 &   0.9794  &  0.9896  &  0.9923 &   0.9917   & 0.9919  &  0.9923   & 0.9923  &  0.9938\\
\hline
 0.6062 &   0.9805 &  0.9314  &  0.9462  &  0.9706  &  0.9836  &  0.9853   & 0.9847  &  0.9848   & 0.9862  &  0.9871\\
 \hline
    0.5368 &   0.8422  &  0.8550  &  0.8911   & 0.9394 &   0.9599  &  0.9684  &  0.9741 &   0.9727 &   0.9793  &  0.9776\\
\hline
    0.4713  &  0.6417  &  0.7334   & 0.8186  &  0.8622  &  0.9154  &  0.9448  &  0.9501  &  0.9505  &  0.9665  &  0.9637\\
 \hline
    0.0000  &  0.4184  &  0.5962  &  0.7144  &  0.8053  &  0.8825  &  0.9044   & 0.9210  &  0.9320  &  0.9544   & 0.9534\\
    \hline
    0.4742  &  0.5068   & 0.6224  &  0.7239  &  0.7988  &  0.8726  &  0.8860  &  0.9097   & 0.9263  &  0.9475  &  0.9580\\
    \hline
    0.5410   & 0.6003  &  0.6988  &  0.7585  &  0.8236 &   0.8860  &  0.9128   & 0.9257  &  0.9358  &  0.9452  &  0.9528\\
    \hline
    0.6106  &  0.6826  &  0.7416  &  0.8125  &  0.8725  &  0.9241  &  0.9305  &  0.9375  &  0.9507  &  0.9567  &  0.9604\\
    \hline
    0.6864  &  0.7330  &  0.7979  &  0.8588  &  0.9183  &  0.9473 &   0.9481  &  0.9480  &  0.9559  &  0.9634 &   0.9733\\
    \hline
    0.7462  &  0.8032 &   0.8564   & 0.9138  &  0.9394   & 0.9610  &  0.9678 &   0.9678   & 0.9696  &  0.9755  &  0.9764\\
    \hline
    \hline
\end{tabular}
\end{center}
\end{table}}
In Figures \ref{figinte2}, \ref{figinte3}, and \ref{figinte5} one displays the level sets of the ratii $\inf\Q_r/T(x_0)$
for $r=2,3$ and $5$ respectively.
The closer to white the color, the closer to $1$ the ratio $\inf\Q_r/T(x_0)$.

\begin{figure}
\resizebox{0.67\textwidth}{!}
{\includegraphics{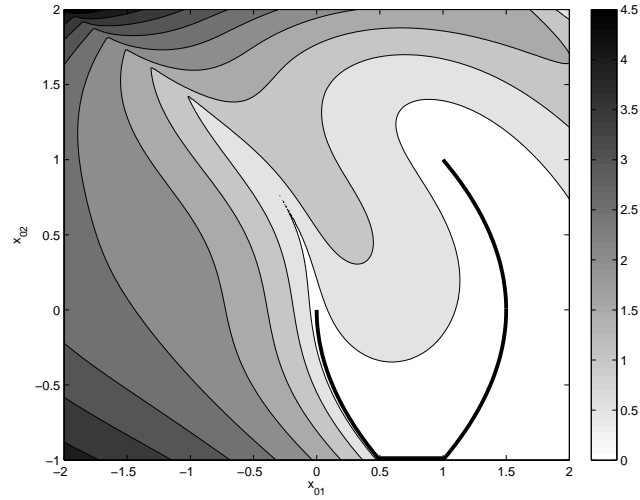}}
\caption{Double integrator: level sets $\Lambda_5(x)-T(x)$ and
optimal trajectory starting from $x_1(0)=x_2(0)=1$.}
\label{figinte6}
\end{figure}

Finally, for this double integrator example we have plotted in Figure \ref{figinte6} the level sets of the function
$\Lambda_5(x)-T(x)$ where $T(x)$ is the known optimal
minimum time to steer the initial state $x$ to $0$; the problem being time-homogeneous, one may take $\Lambda_r\in\RR[x]$ instead of $\RR[t,x]$.
For instance, one may verify that when the initial state is in
the region where the approximation is good, then the whole optimal trajectory also lies in that region.

\subsection{The Brockett integrator}
\label{sec2}  

Consider the so-called \textit{Brockett system} in ${\mathbb{R}}^3$
\begin{equation}\label{brosys}
\begin{array}{rcl}
\dot{x}_1(t)&=&u_1(t),
\\
\dot{x}_2(t)&=&u_2(t),
\\
\dot{x}_3(t)&=&u_1(t)x_2(t)-u_2(t)x_1(t),
\\
\end{array}
\end{equation}
where $x=(x_1,x_2,x_3)$, and the control $\mathbf{u}=(u_1(t),u_2(t))\in\mathcal{U},$
satisfies the constraint
\begin{equation}\label{cont}
u_1(t)^2+u_2(t)^2\leq 1, \quad \forall t\geq0.
\end{equation}
In this case, we have $\X=\RR^3$, $\K=\{(0,0,0)\}$, and $\U$ is the
closed unit ball of $\RR^2$, centered at the origin.

Note that set $\X$ is not compact and so the convergence result
of Theorem \ref{main} may not hold, see the discussion at the beginning of
example \ref{ex1}.  Nevertheless, in the numerical examples,
we have not enforced additional state constraints.
 
\subsubsection{Exact computation}

Let $T(x)$ be the minimum time needed to steer an initial
condition $x\in{\mathbb{R}}^3$ to the origin. We recall the following result of
\cite{BGG} (in fact given for equivalent (reachability) OCP of 
steering the origin to a given point $x$). 

\begin{proposition}
Consider the minimum time OCP for the system (\ref{brosys})
with control constraint (\ref{cont}). The minimum time $T(x)$
needed to steer the origin to a point $x=(x_1,x_2,x_3)\in{\mathbb{R}}^3$
is given by
\begin{equation}\label{formuleT}
T(x_1,x_2,x_3) =
\frac{\theta \sqrt{x_1^2+x_2^2+2\vert x_3\vert}}
{\sqrt{\theta+\sin^2\theta-\sin\theta\cos\theta\:}},
\end{equation}
where $\theta=\theta(x_1,x_2,x_3)$ is the unique solution in
$[0,\pi)$ of
\begin{equation}\label{implic}
\frac{\theta-\sin\theta\cos\theta}{\sin^2\theta} (x_1^2+x_2^2)
=2\vert x_3\vert .
\end{equation}
Moreover, the function $T$ is continuous on ${\mathbb{R}}^3$, and is
analytic outside the line $x_1=x_2=0$.
\end{proposition}

\begin{remark}\label{remensemblecourbes}
Along the line $x_1=x_2=0$, one has
$$T(0,0,x_3)=\sqrt{2\pi\vert x_3\vert}.$$
The singular set of the function $T$, \textit{i.e.}\ the set where $T$
is not
$C^1$, is the line $x_1=x_2=0$ in ${\mathbb{R}}^3$. More precisely, the
gradients $\partial T/\partial x_i$, $i=1,2$, are discontinuous at
every point $(0,0,x_3)$, $x_3\neq 0$.
For the interested reader, the level sets
$\{(x_1,x_2,x_3)\in{\mathbb{R}}^3\ \vert\
T(x_1,x_2,x_3)=r\}$, with $r>0$, are displayed, e.g., in Prieur and Tr\'elat 
\cite{prieur}.
\end{remark}

\subsubsection{Numerical approximation}

Recall that the convergence result of Theorem \ref{main} is 
guaranteed for $\X$ compact only. However, in the present case $\X=\RR^3$ is not compact. One possibility is to take for $\X$ a large ball of $\RR^3$ centered at the origin because
for initial states $x_0$ with norm $\Vert x_0\Vert$
relatively small, the optimal trajectory remains in $\X$.
However, in the numerical experiments presented below, we have chosen to
maintain $\X=\RR^3$, that is, the LMI-relaxation $\Q_r$ does not include
any localizing constraint $M_{r-r(v_j)}(v_j\,z(x))\succeq0$ on the moment sequence $z(x)$.

Recall that $\inf \Q_r\uparrow \min\P$ as $r$ increases, i.e.,
the more moments we consider, the closer to the exact value we get.

In Table \ref{table_example2} we have
displayed the optimal values $\inf\Q_r$ 
for $16$ different values of the
initial state $x(0)=x_0$, in fact, all $16$ combinations of $x_{01}=0$,
$x_{02}=0,1,2,3$, and $x_{03}=0,1,2,3$. So, the entry $(2,3)$ of Table
\ref{table_example2} for the second LMI-relaxation 
is $\inf{\rm Q}_2$ for the initial condition $x_0=(0,1,2)$. At some (few) places in the table, the $^*$ indicates that the SDP solver encountered some numerical problems, which explains why one finds a lower bound $\inf\Q_{r-1}$ slightly higher than 
$\inf\Q_{r}$, when practically equal to the exact value $T^*$.

\begin{table}
\caption{Brockett integrator: LMI-relaxations: $\inf {\rm Q}_r$}
\label{table_example2}
\begin{center}
\begin{tabular}{|| l  | l   | l   | l || }
\hline
\multicolumn{4}{||c||}{first LMI-relaxation: r=1}\\
\hline
    0.0000 &   0.9999 &   1.9999 &   2.9999\\
    \hline
    0.0140 &   1.0017 &   2.0010 &   3.0006\\
    \hline
    0.0243  &  1.0032  &  2.0017  &  3.0024\\
    \hline
    0.0295   & 1.0101   & 2.0034    & 3.0040\\
\hline
\hline
\multicolumn{4}{||c||}{Second LMI-relaxation: r=2}\\
\hline
\hline
0.0000 &   0.9998  &  $1.9997^*$  &  $2.9994^*$ \\
\hline
    0.2012 &   1.1199  &  2.0762 &   3.0453\\
\hline
    0.3738 &   1.2003  &  2.1631 &   3.1304\\
\hline
    0.4946  &  1.3467   & 2.2417  &  3.1943 \\
\hline
\hline
\multicolumn{4}{||c||}{Third LMI-relaxation: r=3}\\
\hline
\hline
0.0000   & 0.9995  &  $1.9987^*$ &   $2.9984^*$ \\
\hline
    0.7665  &  1.3350   & 2.1563 &   3.0530\\
\hline
        1.0826 &   1.7574   & 2.4172   & 3.2036\\
\hline
            1.3804  &  2.0398  &  2.6797 &   3.4077 \\
\hline
\hline
\multicolumn{4}{||c||}{Fourth LMI-relaxation: r=4}\\
\hline
\hline
0.0000  &  0.9992   & 1.9977  &  2.9952\\
\hline
    1.2554 &   1.5925 &   2.1699  &  3.0478\\
    \hline
    1.9962  &  2.1871 &   2.5601  &  3.1977\\
\hline
    2.7006  &  2.7390  &  2.9894 &   3.4254 \\
\hline
\hline
\multicolumn{4}{||c||}{Optimal time $T^*$}\\
\hline
\hline
     0.0000  &  1.0000  &  2.0000   & 3.0000\\
     \hline
    2.5066    &1.7841 &   2.1735  &  3.0547\\
          \hline
    3.5449  &  2.6831 &   2.5819  &  3.2088\\
               \hline
    4.3416   & 3.4328 &   3.0708&    3.4392 \\
  \hline
\end{tabular}
\end{center}
\end{table}

Notice that the upper triangular part (i.e., when both first
coordinates $x_{02},x_{03}$ of the initial condition are away from
zero) displays very good approximations with low order moments.
In addition, the further the coordinates from zero, the best.

For another set of initial conditions
$x_{01}=1$ and $x_{0i}=\{1,2,3\}$ Table \ref{table_example3} displays the results obtained at the LMI-relaxation $\Q_4$.

\begin{table}
\caption{Brockett integrator: $\inf\Q_4$ with $x_{01}=1$}
\label{table_example3}
\begin{center}
\begin{tabular}{||c|c|c||}
\hline
\multicolumn{3}{||c||}{fourth LMI-relaxation: r=4}\\
\hline
1.7979  &  2.3614  &  3.2004 \\
  \hline
 2.3691 &   2.6780  &  3.3341\\
 \hline
  2.8875  &  3.0654   & 3.5337 \\
\hline
\hline
\multicolumn{3}{||c||}{Optimal time $T^*$}\\
\hline
\hline
1.8257 &   2.3636  &  3.2091\\
\hline
    2.5231 &   2.6856 &   3.3426\\
\hline
    3.1895  &  3.1008   & 3.5456 \\
  \hline
\end{tabular}
\end{center}
\end{table}

The regularity property of the minimal-time
function seems to be an important topic of further investigation.

\subsection{Certificate of uncontrollabilty in finite time}
\label{controllability}

Consider the so-called Zermelo problem in $\RR^2$ studied in Bokanowski et al. \cite{boka}
\begin{equation}
\label{zermelo}
\begin{array}{lcl}
\dot x_1(t)&=&1-a\,x_2(t)+v\,\cos\theta\\
\dot x_2(t)&=&v\,\sin\theta\end{array}
\end{equation}
with $a=0.1$. The state $x$ is constrained to remain in the set $\X:=[-6,2]\times [-2,2]\subset\RR^2$, and we also have the control constraints 
$0\leq v\leq 0.44$, as well as $\theta\in [0,2\pi]$. The target $\K$ to reach from an initial state $x_0$ is the ball $B(0,\rho)$ with $\rho:=0.44$. 

With the change of variable $u_1=v\cos \theta$, $u_2=v\sin\theta$,
and $\U:=\{u\::\: u_1^2+u_2^2\leq \rho^2\}$,
this problem is formulated as a minimum time OCP 
with associated hierarchy of LMI-relaxations $(\Q_r)$ defined in (\ref{sdpqrt}). Therefore, if an LMI-relaxation $\Q_r$ at some stage $r$ of the hierarchy is infeasible then the OCP itself is infeasible, i.e., the initial state $x_0$ cannot be steered to the target $\K$ while the trajectory remains in $\X$.

\begin{figure}
\resizebox{0.67\textwidth}{!}
{\includegraphics{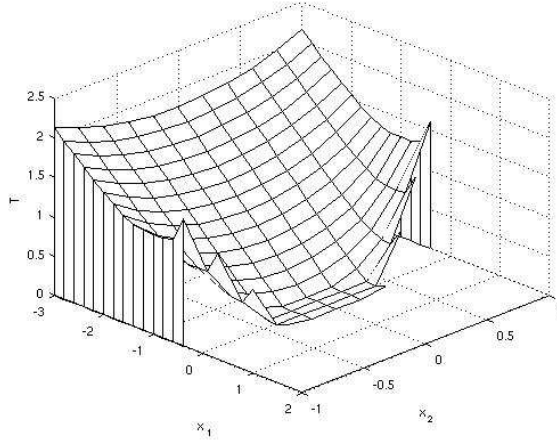}}
\caption{Zermelo problem: uncontrollable states with $\Q_1$}
\label{zermelo1}
\end{figure}

\begin{figure}
\resizebox{0.67\textwidth}{!}
{\includegraphics{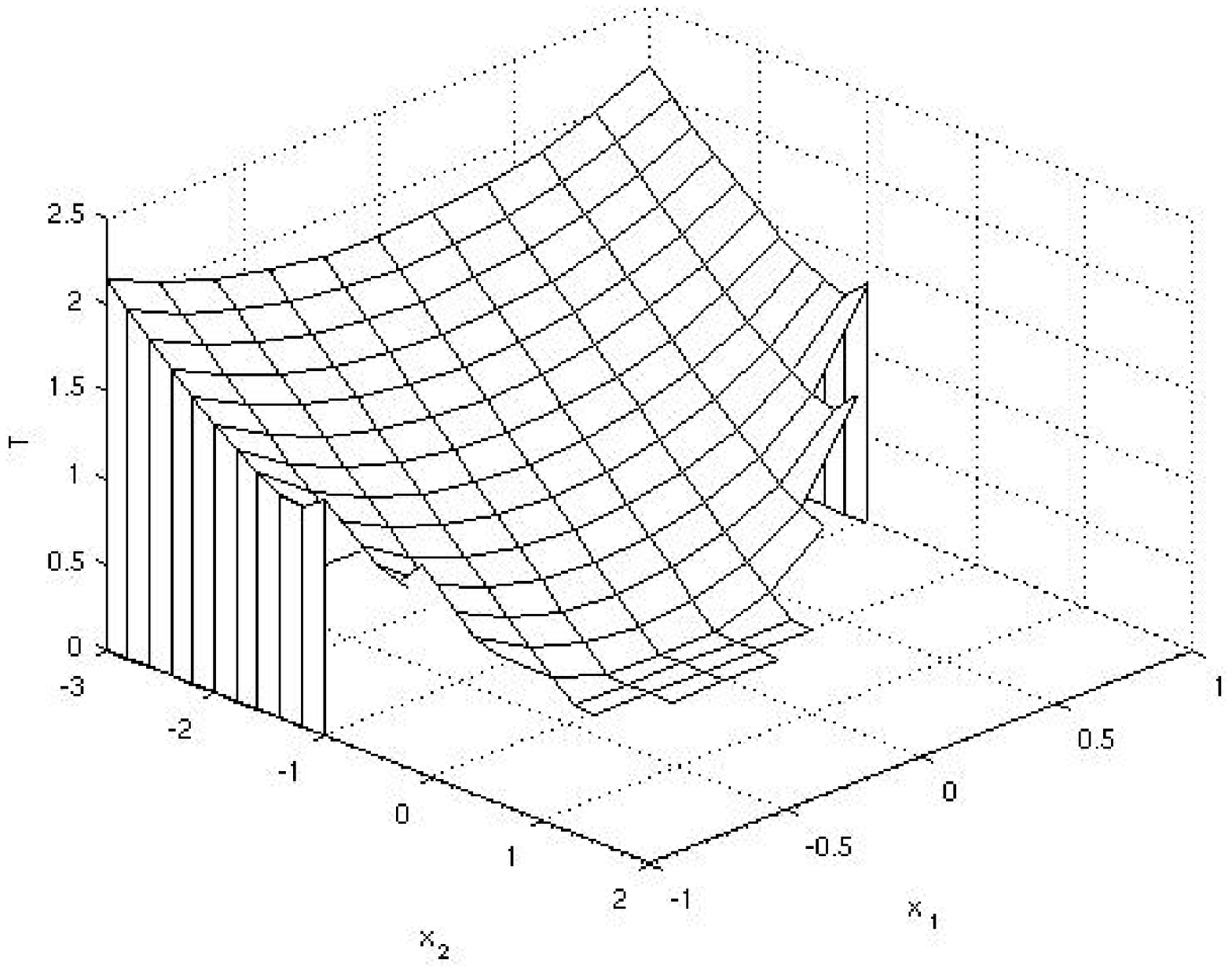}}
\caption{Zermelo problem: uncontrollable states with $\Q_2$}
\label{zermelo2}
\end{figure}

Figures \ref{zermelo1} and \ref{zermelo2} display the uncontrollable initial states $x_0\in\X$ found with the infeasibility of 
the LMI-relaxations $\Q_1$ and $\Q_2$ respectively. 
With $\Q_1$, i.e. by using only second moments,
we already have a very good approximation
of the controllable set displayed in \cite{boka}, and $\Q_2$ brings only a small additional
set of uncontrollable states.

\newpage

\section*{Appendix}
\label{appendix}

\subsection{Proof of Theorem \ref{thlp}}
\label{sec-thlp}
We first prove Item (i).
Consider the linear program $\P$ defined in (\ref{p}), assumed to
be feasible. From the constraint
$\mathcal{L}^*(\mu,\nu)=\delta_{(0,x_0)}$, one has
$$\int_{\K} g(T,x)d\nu - \int_{\s}\left(\frac{\partial g}{\partial t}(t,x) +
  \langle\frac{\partial g}{\partial x}(t,x),f(t,x,u)\rangle\right)d\mu
= g(0,x_0),\quad\forall g\in C_1(\Si).$$
In particular, taking $g(t,x)=1$ and $g(t,x)=T-t$, it follows that
$\mu(\s)=T$ and $\nu(\K)=1$. Hence, for every $(\mu,\nu)$ satisfying
$\mathcal{L}^*(\mu,\nu)=\delta_{(0,x_0)}$,
the pair $(\frac{1}{T}\mu,\nu)$ belongs to the unit ball $B_1$
of $(\m(\s)\times\m(\K))$. From Banach-Alaoglu Theorem, $B_1$ is
compact for the weak $\star$ topology, and even sequentially compact
because $B_1$ is metrizable
(see e.g.\ Hern\'andez-Lerma and Lasserre \cite[Lemma 1.3.2]{ohl}).
Since $\mathcal{L}^*$ is continuous (see Remark
\ref{continuity}), it follows that the set of $(\mu,\nu)$ satisfying
$\mathcal{L}^*(\mu,\nu)=\delta_{(0,x_0)}$ is a
closed subset of $B_1\cap (\m(\s)_+\times\m(\K)_+)$, and thus is
compact. Moreover, since the linear program $\P$ is feasible, this set
is nonempty. Finally, since the linear functional to be minimized is
continuous, $\P$ is solvable.
 
We next prove Item (ii). Consider the set
\[D\,:=\,\{ (\mathcal{L}^*(\mu,\nu),\langle (h,H),(\mu,\nu)\rangle )\ \ 
\vert\ \ (\mu,\nu)\in\m(\s)_+\times\m(\K)_+\}.\]
To prove that $D$ is closed, consider a sequence
$\{(\mu_n,\nu_n)\}_{n\in\N}$ of $\m(\s)_+\times\m(\K)_+$ such that 
\begin{equation}\label{net}
(\mathcal{L}^*(\mu_\alpha,\nu_\alpha),\langle (h,H),(\mu_\alpha,\nu_\alpha)\rangle)\,\to\,(a,b),
\end{equation}
for some $(a,b)\in C_1(\Si)^*\times \RR$. It means that
$\mathcal{L}^*(\mu_n,\nu_n)\to a$,
and $\langle (h,H),(\mu_n,\nu_n)\rangle\to b$.
In particular, taking $\varphi_0:=T-t$ and $\varphi_1=1$,
there must hold
\[\mu_n(\s)\,=\,\langle\varphi_0,\mathcal{L}^*(\mu_n,\nu_n)\rangle
\to \langle\varphi_0,a\rangle,\quad
\nu_n(\K)\,=\,\langle\varphi_1,\mathcal{L}^*(\mu_n,\nu_n)\rangle
\to \langle\varphi_1,a\rangle.\]
Hence, there exist $n_0\in \N$ and a ball $B_r$ of
$\m(\s)\times\m(\K)$, such that $(\mu_n,\nu_n)\in B_r$ for
every $n\geq n_0$. Since $B_r$ is compact, up to a subsequence
$(\mu_n,\nu_n)$ converges weakly to some $(\mu,\nu)\in
\m(\s)_+\times\m(\K)_+$. This fact, combined with
(\ref{net}) and the continuity of $\mathcal{L}^*$, yields
$a=\mathcal{L}^*(\mu,\nu)$, and $b=\langle(h,H),(\mu,\nu)\rangle$.
Therefore, the set $D$ is closed.\\

From Anderson and Nash \cite[Theorems 3.10 and 3.22]{anderson}, it
follows that there is no duality gap between $\P$ and $\P^*$, and
hence, with (i), $\sup\P^*=\min\P$.

Item (iii) follows from Vinter \cite[Theorems 2.1 and 2.3]{vinter},
applied to the mappings
\[F(t,x)\,:=\,f(t,x,U)\,,\quad l(t,x,v)\,:=\,\inf_{u\in
  U}\:\{\:h(t,x,u)\ \vert \ v=f(t,x,u)\:\},\]
 for $(t,x)\in \RR\times \RR^n$.
 $\qed$

\subsection{Proof of Theorem \ref{main}}
\label{proofmain}

First of all, observe that $\Q_r$ has a feasible solution. Indeed, it
suffices to consider the sequences $y$ and $z$ consisting of the
moments (up to order $2r$) of the occupations measures
$\nu^\mathbf{u}$ and $\mu^\mathbf{u}$ associated with an admissible
control $\mathbf{u}\in\mathcal{U}$ of the OCP
(\ref{system})-(\ref{crit1}).

Next, observe that, for every couple $(y,z)$ satisfying all
constraints of $\Q_r$, there must holds $y_0=1$ and $z_0=1$.
Indeed, it suffices to choose $g(t,x)=1$ and $g(t,x)=1-t$ (or $t$) in the
constraint
\[L_y(g_1) -L_z(\partial g/\partial t+\langle \nabla_x
g,f\rangle)=g(0,x_0).\]

We next prove that, for $r$ sufficiently large, one has
$\vert z(x)_\alpha\vert\leq 1$, $\vert z(u)_\beta\vert\leq 1$,
$\vert z(t)_k\vert\leq 1$, for every $k$, and $\vert
y_\alpha\vert\leq1$. We only provide the details of the proof for
$z(x)$, the arguments being similar for $z(u)$, $z(t)$ and $y$.

Let $\Sigma^2\subset\RR[x]$ be the space of sums of squares (s.o.s.)
polynomials, and let $Q\subset\RR[x]$ be the {\it quadratic modulus}
generated by the polynomials $v_j\in\RR[x]$ that define $\X$, i.e.,
\[Q\,:=\,\{\:\sigma\in\RR[x]\quad\vert\quad \sigma\,=\,\sigma_0+\sum_{j\in J}\sigma_j\,v_j\quad \mbox{with }\sigma_j\in\Sigma^2,\quad \forall\,j\in \{0\}\cup J\}.\]
In addition, let  $Q(t)\subset Q$ be the set of elements $\sigma$ of
$Q$ which have a representation $\sigma_0+\sum_{j\in J}\sigma_j\,v_j$
for some s.o.s. family $\{\sigma_j\}\subset\Sigma^2$ with ${\rm
  deg}\,\sigma_0\leq 2t$ and ${\rm deg}\,\sigma_jv_j\leq 2t$ for every
$j\in J$.

Let $r\in\N$ be fixed. Since $\X\subset [-1,1]^n$, there holds
$1\pm x^\alpha >0$ on $\X$, for every $\alpha\in\N^n$ with
$\vert\alpha\vert\leq 2r$.
Therefore, since $\X$ satisfies Putinar' condition (see Definition
\ref{def1}), the polynomial  $x\mapsto 1\pm x^\alpha$ belongs to $Q$
(see Putinar \cite{putinar}). Moreover, there exists $l(r)$ such that
$x\mapsto 1\pm x^\alpha\in Q(l(r))$ for every $\vert\alpha\vert\leq 2r$.
Of course, $x\mapsto 1\pm x^\alpha\in Q(l)$ for every
$\vert\alpha\vert\leq 2r$, whenever  $l\geq l(r)$.

For every feasible solution $z$ of $\Q_{l(r)}$ one has
\[\vert z(x)_\alpha\vert\,=\,\vert\: L_{z}(x^\alpha)\:\vert\leq
z_0=1,\qquad \vert\alpha\vert \leq 2r.\]
This follows from $z_0=1$, $M_{l(r)}(z)\succeq0$ and
$M_{l(r)-r(v_j)}(v_j\,z(x))\succeq0$, which implies
\[z_0\pm z(x)_\alpha\,=\,L_{z}(1\pm x^\alpha)\,=\,L_{z}(\sigma_0)+ 
\sum_{j=1}^mL_{z(x)}(\sigma_j\,v_j)\geq0.\]
With similar arguments, one redefines $l(r)$ so that, for every couple
$(y,z)$ satisfying the contraints of $\Q_{l(r)}$, one has
\[0\leq z_k(t)\leq 1\quad\mbox{and}\quad\vert z(x)_\alpha\vert,\:\vert z(u)_\beta\vert,\:\vert y_\alpha\vert\,\leq\,1,\qquad \forall \,k,\:\vert\alpha\vert,\:\vert\beta\vert \leq 2r.\]

From this, and with $l(r):=2l(r)$, we immediately deduce that
$\vert z_\gamma\vert\leq 1$ whenever $\vert\gamma\vert\leq 2r$.
In particular, $L_y(H)+L_z(h)\,\geq\,-\sum_\beta \vert H_\beta\vert 
-\sum_\gamma \vert h_\gamma\vert$, 
which proves that $\inf\Q_{l(r)}>-\infty$, and so $\inf\Q_r>-\infty$
for $r$ sufficiently large.\\
Let $\rho:=\inf\P=\min\P$ (by Theorem \ref{thlp}), let $r\geq l(r_0)$, and let $(z^r,y^{r})$ be a nearly optimal solution of 
$\Q_r$ with value
\begin{equation}
\label{value}
\inf\Q_{r}\leq L_{z^r}(h)+L_{y^{r}}(H)\leq\inf\Q_{r}+\frac{1}{r}\quad\left(\leq\rho+\frac{1}{r}\right).
\end{equation}
Complete the finite vectors $y^r$ and $z^r$ with zeros to make them
infinite sequences. Since for arbitrary $s\in\N$ one has
$\vert y^r_\alpha\vert,\,\vert z^r_\gamma\vert\leq 1$ whenever
$\vert\alpha\vert,\vert\gamma\vert\leq 2s$, provided $r$ is sufficiently large,
by a standard diagonal argument, there exists a subsequence $\{r_k\}$ and two infinite sequences 
$y$ and $z$, with $\vert y_\alpha\vert\leq1$ and $\vert z_\gamma\vert\leq1$, for all 
$\alpha,\gamma$, and such that
\begin{equation}
\label{convergence}
\lim_{k\to\infty}\,y_\alpha^{r_k}\,=\,y_\alpha\quad\forall \alpha\in\N^n;\quad
\lim_{k\to\infty}\,z_\gamma^{r_k}\,=\,z_\gamma\quad\forall \gamma\in\N\times\N^n\times\N^m.
\end{equation}
Next, let $r$ be fixed arbitrarily. Observe that 
$M_{r_k}(y^{r_k})\succeq0$ implies $M_{r}(y^{r_k})\succeq0$
whenever $r_k\geq r$, and similarly,
$M_{r}(z^{r_k})\succeq0$. Therefore, from (\ref{convergence}) and 
$M_{r}(y^{r_k})\succeq0$, we deduce that $M_r(y)\succeq0$, and similarly $M_r(z)\succeq0$. Since this holds for arbitrary $r$, and 
$\vert y_\alpha\vert,\,\vert z_\gamma\vert\leq 1$ for all
$\alpha,\gamma$, one infers from Proposition \ref{prop1}
that $y$ and $z$ are moment sequences of two measures
$\nu$ and $\mu$ with support contained in $[-1,1]^n$ and
$[0,1]\times [-1,1]^n\times [-1,1]^m$ respectively. In addition, 
from the equalities $y^{r_k}_0=1$ and $z^{r_k}_0=1$ 
for every $k$, it follows that $\nu$ and $\mu$ are probability measures
on $[-1,1]^n$, and $[0,1]\times [-1,1]^n\times [-1,1]^m$.
 
Next, let $(t,\alpha)\in\N\times\N^n$ be fixed, arbitrarily. From
\[L_{y^{r_k}}(g_1) -g(0,x_0) -L_{z^{r_k}}(\partial g/\partial t+\langle \nabla_x g,f\rangle)=0,\quad\mbox{with }g=(t^px^\alpha),\]
and the convergence (\ref{convergence}), we obtain
\[L_{y}(g_1) -g(0,x_0) -L_{z}(\partial g/\partial t+\langle \nabla_x g,f\rangle)=0,\quad\mbox{with }g=(t^px^\alpha), \]
that is, $\langle \mathcal{L}g ,(\mu,\nu)\rangle\,=\,\langle g,\delta_{(0,x_0)}\rangle$.
Since $(t,\alpha)\in\N\times\N^n$ is arbitrary, we have
\[\langle g,\mathcal{L}^*(\mu,\nu)\rangle\,=\,
\langle \mathcal{L}\,g,(\mu,\nu)\rangle\,=\,\langle g,\delta_{(0,x_0)}\rangle\quad
\forall \,g\in\RR[t,x],\]
which implies that $\mathcal{L}^*(\mu,\nu)=\delta_{(0,x_0)}$.

Let $z(x)$, $z(u)$ and $z(t)$ denote the moment vectors of
the marginals of $\mu$ with respect to the variables $x\in\RR^n$,
$u\in\RR^m$ and $t\in\RR$, respectively, i.e.,
\[z(x)^\alpha=\int x^\alpha\,\mu(d(t,x,u))\quad \forall\,\alpha\in\N^n,\quad
z(u)^\beta=\int u^\beta\,\mu(d(t,x,u))\quad \forall\beta\in\N^m,\]
and $z(t)^k=\int t^k\,\mu(d(t,x,u))$ for every $k\in\N$.

With $r$ fixed arbitrarily, and using again (\ref{convergence}),
we also have
$M_r(\theta_jy)\succeq0$ for every $j\in J_T$, and
\[ M_{r}(v_j\,z(x))\succeq0\quad\forall \,j\in J,\quad M_{r}(w_k\,z(u))\succeq0\quad\forall k\in W,
\quad M_{r}(t(1-t)\,z(t))\succeq0.\]

Since $\X$, $\K$ and $\U$ satisfy Putinar's condition (see Definition
\ref{def1}), from Theorem \ref{putinarthm} (Putinar's Positivstellensatz), $y$ 
is the moment sequence of a probability measure
$\nu$ supported on $\K\subset [-1,1]^n$. Similarly,
$z(x)$ is the moment sequence of a measure
$\mu^x$ supported on $\X\subset [-1,1]^n$,
$z(u)$ is the moment sequence of a measure
$\mu^u$ supported on $\U\subset [-1,1]^m$, and
$z(t)$ is the moment sequence of a measure
$\mu^t$ supported on $[0,1]$. Since measures on compact sets are
moment determinate, it follows that $\mu^x$, $\mu^u$, and $\mu^t$ are
the marginals of $\mu$ with respect to $x$, $u$ and $t$ respectively.
Therefore, $\mu$ has its support contained in $\s$, and from
$\mathcal{L}^*(\mu,\nu)=\delta_{(0,x_0)}$ it follows that $(\mu,\nu)$
satisfies all constraints of the problem $\P$.

Moreover, one has
\begin{eqnarray*}
\lim_{k\to\infty}\,\inf\Q_{r_k}&=&\lim_{k\to\infty}\,L_{z^{r_k}}(h)+L_{y^{r_k}}(H)
\qquad\mbox{(by (\ref{value}))}\\
&=&L_{z}(h)+L_{y}(H)\qquad\mbox{(by (\ref{convergence}))}\\
&=&\int h\,d\mu+\int H\,d\nu\quad\leq\rho=\min\P.
\end{eqnarray*}
Hence, $(\mu,\nu)$ is an optimal solution of $\P$, and
$\min\Q_r\uparrow\min\P$ (the sequence is monotone nondecreasing).
Item (i) is proved.

Item (ii) follows from Theorem \ref{thlp} (iii).
 $\qed$

\subsection{Proof of Theorem \ref{mainbis}}
\label{proofmainbis}
It suffices to prove that $v_p\to v$ as $p\to +\infty$. For every
integer $p$, $v_p=\min\P_p$ is attained for a couple of measures
$(\mu_p,\nu_p)$. As in the proof of Theorem \ref{thlp}, the sequence
$\{(\mu_p,\nu_p)\}_{p\in\N}$ is bounded in $\m(\s)_+\times\m(\K)_+$,
and hence, up to a subsequence, it converges to an element $(\mu,\nu)$
of this space for the weak $\star$ topology.

On the one hand, by definition, 
$\mathcal{L}_p^*(\mu_p,\nu_p)=\delta_{(0,x_0)}$ for every $p$.
On the other, $\mathcal{L}_p^*$ tends
strongly to $\mathcal{L}^*$, and so
$\mathcal{L}^*(\mu,\nu)=\delta_{(0,x_0)}.$
Moreover, since $(h_p,H_p)$ tends strongly to $(h,H)$ in
$C_1(\s)\times C_1(\K)$, one has
$$v_p=\langle(\mu_p,\nu_p),(h_p,H_p)\rangle \ \longrightarrow
\langle(\mu,\nu),(h,H)\rangle,$$
and so $v \leq \langle(\mu,\nu),(h,H)\rangle$.
We next prove that $v = \langle(\mu,\nu),(h,H)\rangle$.

Since $(\mu_p,\nu_p)$ is an optimal solution of $\P_p$,
$$\langle(\mu_p,\nu_p),(h_p,H_p)\rangle \ \leq\
\langle(\bar\mu,\bar\nu),(h_p,H_p),\quad\forall(\bar\mu,\bar\nu)\
\vert\ \mathcal{L}_p^*(\bar\mu,\bar\nu)=\delta_{(0,x_0)}.$$
Hence, passing to the limit,
$$\langle(\mu,\nu),(h,H)\rangle \ \leq\
\langle(\bar\mu,\bar\nu),(h,H),\quad\forall(\bar\mu,\bar\nu)\
\vert\ \mathcal{L}^*(\bar\mu,\bar\nu)=\delta_{(0,x_0)},$$
and so, $(\mu,\nu)$ is an optimal solution of $\P$, i.e.,
$v = \langle(\mu,\nu),(h,H)\rangle$.
$\qed$

\section*{Acknowledgments}

The research of Didier Henrion was partly supported by
Project No.~102/06/0652 of the Grant Agency of the Czech Republic
and Research Program No.~MSM6840770038 of the Ministry
of Education of the Czech Republic.


\end{document}